\documentclass[11pt]{amsart}
\theoremstyle{plain}
\newtheorem{thm}{Theorem}[section]
\newtheorem{theorem}[thm]{Theorem}

\newtheorem{lemma}[thm]{Lemma}
\newtheorem{corollary}[thm]{Corollary}
\newtheorem{proposition}[thm]{Proposition}
\theoremstyle{definition}
\newtheorem{remark}[thm]{Remark}

\newtheorem{notation}[thm]{Notation}

\newtheorem{definition}[thm]{Definition}

\newtheorem{example}[thm]{Example}

\numberwithin{equation}{section}
\newcommand{\Fr}{{\mathbf F}}

\newcommand{\ad}{{\rm ad}}

\newcommand{\p}{\partial}
\newcommand{\ep}{\varepsilon}
\newcommand{\rd}{{\rm d}}
\newcommand{\tH}{\textsf{H}}

\newcommand{\sA}{{\mathcal A}}

\newcommand{\sC}{{\mathcal C}}

\newcommand{\sE}{{\mathcal E}}
\newcommand{\sF}{{\mathcal F}}
\newcommand{\sG}{{\mathcal G}}
\newcommand{\sH}{{\mathcal H}}

\newcommand{\sJ}{{\mathcal J}}
\newcommand{\sK}{{\mathcal K}}

\newcommand{\sO}{{\mathcal O}}
\newcommand{\sP}{{\mathcal P}}

\newcommand{\sT}{{\mathcal T}}

\newcommand{\C}{{\mathbb C}}

\newcommand{\BK}{{\mathbb K}}

\newcommand{\BP}{{\mathbb P}}

\newcommand{\Q}{{\mathbb Q}}

\newcommand{\End}{{\rm End}}

\newcommand{\fg}{{\mathfrak g}}
\newcommand{\fb}{{\mathfrak b}}
\newcommand{\fsl}{{\mathfrak s}{\mathfrak l}}
\newcommand{\fgl}{{\mathfrak g}{\mathfrak l}}
\newcommand{\ff}{{\mathfrak f}}

\newcommand{\fe}{{\mathfrak e}}
\newcommand{\ft}{{\mathfrak t}}
\newcommand{\fl}{{\mathfrak l}}

\newcommand{\fu}{{\mathfrak u}}
\newcommand{\fh}{{\mathfrak h}}

\newcommand{\fso}{{\mathfrak s}{\mathfrak o}}

\newcommand{\aut}{{\mathfrak a}{\mathfrak u}{\mathfrak t}}

\newcommand\Aut{\rm Aut}

\newcommand{\wt}{{\rm wt}}

\def\Sym{\mathop{\rm Sym}\nolimits}

\def\Hom{\mathop{\rm Hom}\nolimits}

\title[Conic connections and principal connections]{Characteristic conic connections  and\\ torsion-free principal connections}

\author{Jun-Muk Hwang and Qifeng Li}

\thanks{Jun-Muk Hwang was supported by the Institute for Basic Science (IBS-R032-D1). Qifeng Li was supported by the NSFC grant No. 12201348.}
\begin{document}

\begin{abstract}
We study the relation between torsion tensors of principal connections on G-structures and characteristic conic connections on associated cone structures. We formulate sufficient conditions under which the existence of a characteristic conic connection implies the existence of a torsion-free principal connection. We verify these conditions for adjoint varieties of simple Lie algebras, excluding those of type $\textsf{A}_{\ell \neq 2}$ or $\textsf{C}_{\ell}$. As an application, we give a complete classification of the germs of minimal rational curves whose VMRT at a general point is such an adjoint variety: nontrivial ones come from lines on  hyperplane sections of  certain Grassmannians or minimal rational curves on wonderful group compactifications.
\end{abstract}

\keywords{G-structures, adjoint varieties, minimal rational curves, varieties of minimal rational tangents}

\maketitle

\medskip
MSC2020: 53C10, 53C35, 14M27

\bigskip
{\bf Convention} \begin{itemize} \item[1.] We work in the holomorphic category: all varieties are complex analytic and all maps are holomorphic. Open sets refer to Euclidean topology and open sets in Zariski topology are called Zariski-open sets. A general point of a complex analytic set $X$ means a point in a  Zariski-open subset of $X$, i.e., the complement of a closed analytic subset in $X$.
\item[2.]
We use the following notations regarding projective space.  For a vector space $V$, its projectivization $\BP V$ is the set of 1-dimensional subspaces of $V$ with the natural quotient map $q: V \setminus 0 \to \BP V$.  For a subvariety $Z \subset \BP V$, we write $$ \check{Z} :=  q^{-1}(Z) \   \subset V \setminus 0 \mbox{ and } \widehat{Z} :=  \check{Z} \cup 0  \ \subset V. $$ For a nonzero vector $w \in V$, the corresponding  point of $\BP V$ is denoted by $[w]$ so that $\C w = \widehat{[w]}$.   If $Z$ is a submanifold and $\alpha \in Z$,  the affine tangent space $T_{z} \check{Z} \subset V$  at a point $z \in \check{\alpha}$ is sometimes written as $T_z \widehat{Z}$ or $T_{\alpha} \widehat{Z} \subset V$.
 For a nonnegative integer $k$, the $k$-th tensor product of the hyperplane line bundle (resp. tautological line bundle) on $\BP V$ is denoted by $\sO(k)$ (resp. $\sO(-k)$) and we often write the restriction $\sO(k)|_Z$ as $\sO(k)$ by abuse of notation. The tensor product of a vector bundle $E$ on $Z$ and $\sO(k)$ is abbreviated to $E(k)$.
For a vector bundle $F$ on a complex manifold $M$, its projectivization $\BP F$ is the fiber bundle the fiber of which at $x \in M$ is the projectivization $\BP F_x$ of the vector space $F_x$.

\item[3.] For a  holomorphic map $f: X \to Y$ between complex manifolds,
its differential at $x\in X$ is denoted by ${\rm d}_x f : T_x X \to T_{f(x)} Y$. When $f$ is submersive, namely, when ${\rm d}_x f$ is surjective for all $x\in X$,  the vector subbundle ${\rm Ker}({\rm d} f) \subset T X$ is denoted by $T^f$.
\item[4.] For a principal $G$-bundle $\sP \to M$ on a complex manifold $M$ and a left action of $G$ on an algebraic variety $A$, the fiber bundle associated with $\sP$ with a standard fiber $A$ is denoted by $\sP \times_G A$. It is the quotient of $\sP \times A$ by the left $G$-action $(w, a) \mapsto (w \cdot g^{-1}, g \cdot a)$ for all $w \in \sP, a \in A$ and $g \in G$. If $A$ is a vector space and $G$ acts linearly on $A$, then $\sP \times_G A$ is a vector bundle on $M$.   \item[5.] We denote by $\Delta$ an open neighborhood of the origin in $\C$. \end{itemize}

\section{Introduction}

A cone structure on a complex manifold $M$ is a closed submanifold $\sC \subset \BP TM$ of the projectivized tangent bundle of $M$ such that the projection $ \phi: \sC \to M$ is a submersion. There are  two major sources of (holomorphic) cone structures: twistor theory in differential geometry/mathematical physics (for example, see  Chapter 1, Section 6 of \cite{Ma}) and the theory of minimal rational curves in algebraic geometry (for example, see Section 4 of \cite{Hw10}). An isotrivial cone structure is  a cone structure $\sC \subset \BP TM$ such that all fibers of $\phi: \sC \to M$ are isomorphic to each other  as projective submanifolds.  More precisely, we call it a $Z$-isotrivial cone structure for a projective submanifold $Z \subset \BP V, \dim V = \dim M$, if  $\phi^{-1}(x) \subset \BP T_x M$ is projectively isomorphic to $Z \subset \BP V$ for any $x \in M$.  A $Z$-isotrivial cone structure is naturally associated to a G-structure with the structure group $\Aut(\widehat{Z}) \subset {\rm GL}(V)$, the linear automorphism group of $Z \subset \BP V$ (see Definition \ref{d.isotrivial}).  Conversely,  for a linear group $G \subset {\rm GL}(V)$ and a projective submanifold $Z \subset \BP^{n-1}$ preserved by the $G$-action, a G-structure on $M$ with the structure group $G$ induces  a $Z$-isotrivial cone structure $\sC \subset \BP TM.$  Thus there exists a natural one-to-one correspondence between $Z$-isotrivial cone structures and G-structures with the structure group $\Aut(\widehat{Z}) \subset {\rm GL}(V)$.

A classical example is the cone of null vectors in   holomorphic conformal geometry studied in \cite{LB}.  The null vectors of a local holomorphic metric tensor on $M$ determine a fiber subbundle $\sC \subset \BP TM$ whose fiber $\sC_x \subset \BP T_xM$ is a nonsingular quadric hypersurface for each $x \in M$. Conversely, any smooth hypersurface $\sC \subset \BP TM$ whose fibers $\sC_x \subset \BP TM$ are smooth quadric hypersurfaces can be realized as the cone of null vectors  of a conformal structure on $M$.

A standard method to study G-structures is the theory of principal connections (see \cite{KN}, \cite{Sc}, \cite{St}). Among invariants of a principal connection,  its torsion tensor is fundamental.  In fact, a basic question one can ask about a G-structure is whether it admits a  torsion-free principal connection (this is equivalent to asking whether the $G$-structure is 1-flat as explained in Section 1 of \cite{Br}). On the other hand, for a cone structure $\sC \subset \BP TM$, we have the notion of a conic connection (Chapter 1, Section 6, Definition 3 of \cite{Ma}), which specifies a collection of distinguished curves on $M$ in the direction of $\sC$. Of particular interest are characteristic conic connections, which arise naturally  in the theory of minimal rational curves in algebraic geometry (see Section \ref{s.vmrt}), but also have been used implicitly in twistor theory. For example, the null geodesics of a  Levi-Civita connection of a holomorphic conformal structure give a characteristic conic connection on the cone structure of null vectors. More generally, regular normal  parabolic geometries of certain types can be regarded as characteristic conic connections of some isotrivial cone structures (see  Theorem 0.2 of \cite{HN}).

Our main goal is to study the correspondence between conic connections on isotrivial cone structures and principal connections on the associated G-structures, with special regard to characteristic conic connections and torsion-free principal connections:

$$\begin{array}{ccc} \mbox{Isotrivial cone structures}
& \stackrel{(\textsf{1})}{\Longleftrightarrow} & \mbox{G-structures} \\ \mbox{Conic connections} & \stackrel{(\textsf{2})}{\Longleftrightarrow} & \mbox{Principal connections}  \\
\mbox{ Characteristic conic conn.} & \stackrel{(\textsf{3})}{\Longleftrightarrow} & \mbox{Torsion-free principal conn.}
 \end{array} $$

\bigskip
We have already mentioned the natural one-to-one correspondence (\textsf{1}).
In the second correspondence (\textsf{2}),  that a principal connection on a G-structure gives rise to a conic connection on the associated cone structure is immediate: the conic connection comes from  geodesics of the principal connection.
The reverse direction in (\textsf{2}) is not automatic, but  can be obtained in many interesting cases by  the following result. See Theorem \ref{t.CtoP} for a more precise statement.

\begin{theorem}\label{t.tensor1}
Let $Z \subset \BP V$ be a closed submanifold satisfying \begin{itemize}
\item[\textsf{(i)}] $H^0(Z, \sO(1)) = V^*$ and
\item[\textsf{(ii)}] the natural homomorphism $\aut(\widehat{Z}) \otimes H^0(Z, \sO(1)) \to H^0(Z, TZ\otimes \sO(1))$ is surjective. \end{itemize}
Then any conic connection on a $Z$-isotrivial cone structure   arises locally as the conic connection given by geodesics of a principal connection on the associated G-structure. \end{theorem}

The third correspondence (\textsf{3}) is much more subtle. Our main result regarding  (\textsf{3}) is the following.  See Theorem \ref{t.Xi} for a more precise statement.

\begin{theorem}\label{t.torsion1}
For a  submanifold $Z \subset \BP V$,  let
 $\Xi_Z \subset \Hom(\wedge^2 V, V)$ be the subspace consisting of homomorphisms $\sigma: \wedge^2 V \to V$ satisfying $\sigma(z, T_z \widehat{Z}) \subset T_z \widehat{Z}$ for any $ z \in \check{Z}$.
Let $\sC \subset \BP TM$ be a $Z$-isotrivial cone structure and let $\sH$ be a principal connection  on the associated G-structure on $M$. Then the conic connection induced by $\sH$ is a characteristic conic connection if and only if the torsion tensor of $\sH$ has values in $\Xi_Z \subset \Hom(\wedge^2 V, V).$ \end{theorem}

Because $0$ is an element of $\Xi_Z$, Theorem \ref{t.torsion1} implies that the  conic connection associated to a torsion-free principal connection is a characteristic conic connection, establishing one direction of the correspondence (\textsf{3}). The reverse direction of (\textsf{3}) usually does not hold.  We obtain, however, the following result under some additional conditions on $Z \subset \BP V$. See Theorem \ref{t.Spencer} for a more precise statement.

\begin{theorem}\label{t.torsionfree1}
Assume that $Z \subset \BP V$ satisfies the two conditions \textsf{(i)} and \textsf{(ii)} of Theorem \ref{t.tensor1} and also \begin{itemize} \item[\textsf{(iii)}] the Spencer homomorphism $\p: \Hom(V, \aut(\widehat{Z})) \to \Hom(\wedge^2 V, V)$ defined by $$\p h (u, v) = h(u) \cdot v - h(v) \cdot u \mbox{ for all } u, v \in V$$ satisfies $\Xi_Z \subset {\rm Im}(\p)$. \end{itemize}
Then for any $Z$-isotrivial cone structure equipped with a characteristic conic connection,
 the associated G-structure  admits locally  a torsion-free principal connection. \end{theorem}

Let us  remark here that although the correspondence (\textsf{1}) is one-to-one, the correspondences (\textsf{2}) and (\textsf{3}), when they hold as in  Theorem \ref{t.torsionfree1}, are far from one-to-one. In fact,  when $Z \subset \BP V$ is a not a linear subspace,  a characteristic conic connection on a $Z$-isotrivial cone structure is unique if it exists (see Theorem 3.1.4 of \cite{HM99}). It is well-known from the classical theory of G-structures that there  exist (locally) many different torsion-free principal connections when the Spencer homomorphism is not injective. Thus (\textsf{3}) is not one-to-one in such a case.  Even when the Spencer homomorphism is injective,  two principal connections whose difference (in the sense of Proposition \ref{p.beta} below)  belongs to $\partial^{-1}(\Xi_Z)$ would  give the same characteristic conic connection by Theorem \ref{t.torsion1}. Thus the correspondence (\textsf{2}) is not one-to-one.

Combining Theorem \ref{t.torsionfree1} with the classification result \cite{MS} of Merkulov and Schwachh\"ofer on torsion-free irreducible G-structures, we obtain the following result. See Corollaries \ref{c.MS} and \ref{c.homogeneous} for a more precise statement and see Definition \ref{d.homogeneous} for the definition of a locally symmetric/homogeneous isotrivial cone structure.

\begin{theorem}\label{t.MS1}
Assume that  $Z \subset \BP V$ is a highest weight variety, namely, the action of $\Aut(\widehat{Z}) \subset {\rm GL}(V)$ on $V$ is an irreducible representation and $Z$ is the unique closed orbit in $\BP V$.
 Assume furthermore that  $Z \subset \BP V$ satisfies \textsf{(i)}, \textsf{(ii)}, \textsf{(iii)} of Theorem \ref{t.torsionfree1} and is not one of the projective submanifolds listed in Propositions \ref{p.KN} and  \ref{p.Legendre}.
Then a $Z$-isotrivial cone structure equipped with a characteristic conic connection is locally symmetric, hence locally homogeneous. \end{theorem}

When  $Z \subset \BP V$ is a highest weight variety,  the condition \textsf{(i)} is always satisfied and  the condition \textsf{(ii)} is often easy to verify.  The condition \textsf{(iii)} of Theorem \ref{t.torsionfree1}  is the key for the strong conclusion of Theorem \ref{t.MS1} and seems to hold only in special cases.

As an explicit example of Theorem \ref{t.MS1}, we prove the following.

\begin{theorem}\label{t.adjoint1}
Let $\fg$ be a complex simple Lie algebra of type different from $\textsf{A}_{\ell\neq 2}$ or $\textsf{C}_{\ell}$ (thus it is of type $\textsf{A}_2, \textsf{B}_{\ell \geq 3}, \textsf{D}_{\ell \geq 4}, \textsf{E}_6, \textsf{E}_7, \textsf{E}_8,  \textsf{F}_4$ or $\textsf{G}_2$). Let $Y\subset \BP \fg$ be the adjoint variety of $\fg$, namely, the highest weight variety of the adjoint representation of $\fg$. Then \begin{itemize} \item[(a)] the natural homomorphism $\aut(\widehat{Y}) \otimes H^0(Y, \sO(1)) \to H^0(Y, TY \otimes \sO(1))$ is surjective; and
\item[(b)] the Spencer homomorphism $\p: \Hom(\fg, \aut(\widehat{Y})) \to \Hom(\wedge^2 \fg,\fg)$ satisfies $\Xi_Y \subset {\rm Im}(\p)$. \end{itemize} Consequently, a $Y$-isotrivial cone structure with a characteristic conic connection is locally symmetric.  \end{theorem}

The proof of (a) in  Theorem \ref{t.adjoint1} is relatively simple, but the proof of (b) is rather involved, using  the natural holomorphic contact structure on $Y$ and the geometry of lines on $Y$ relative to this contact structure.

Theorem \ref{t.adjoint1} does not give a complete classification of $Y$-isotrivial cone structures admitting characteristic conic connections. To do that, we need to classify locally symmetric connections on the G-structure  with the structure group $\Aut(\widehat{Y})$.
By the standard theory of symmetric spaces, the classification follows from the following algebraic result.

\begin{theorem}\label{t.Bianchi}
Let $\fg$ be  a complex simple Lie algebra of type different from $\textsf{A}_{\ell\neq 2}$ or $\textsf{C}_{\ell}$ (as in Theorem \ref{t.adjoint1}). If $f \in \Hom(\wedge^2 \fg, \fg)$ and $\sigma \in \wedge^2 \fg^*$ satisfy \begin{eqnarray}\label{e.Bianch} 0 & = &   [f(x,y), z] + [f(y,z), x] + [f(z,x),y] \\ & & + \sigma(x,y) \, z + \sigma(y,z) \, x + \sigma (z,x) \, y \nonumber \end{eqnarray}  for all $x ,y, z \in \fg$, then $\sigma =0$ and $f$ is a scalar multiple of the Lie bracket $ [\, ,\,]: \wedge^2 \fg \to \fg.$ \end{theorem}

The condition (\ref{e.Bianch}) is the Bianchi identity for the curvature tensors of principal connections (see Proposition \ref{p.BK}) and Theorem \ref{t.Bianchi} determines  the space
 of formal curvature maps (see Definition \ref{d.Bianchi}) of the Lie algebra $\fg \oplus \C {\rm Id}_{\fg}$.  The classification of affine symmetric spaces with irreducible holonomies is a classical result.  But to our knowledge,  the classification of affine symmetric spaces subordinate to   irreducible G-structures, with potentially reducible holonomies, is an open problem. Theorem \ref{t.Bianchi} settles this problem for the Lie algebra $\fg \oplus \C {\rm Id}_{\fg}$ for a simple Lie algebra $\fg$ of type different from $\textsf{A}_{\ell \neq 2}$ or $\textsf{C}_{\ell}$. Our Theorem \ref{t.MS1} suggests that it is worth studying  this question for other Lie algebras, especially those appearing in the classification of  symmetric spaces with irreducible holonomies.

 Theorem \ref{t.adjoint1} and Theorem \ref{t.Bianchi} have  the following application in the study of minimal rational curves in algebraic geometry. See Theorem \ref{t.vmrt} for a more precise statement.

\begin{theorem}\label{t.vmrt1}
Let $\fg$ be  a complex simple Lie algebra of type different from  $\textsf{A}_{\ell\neq 2}$ or $\textsf{C}_{\ell}$ and let $Y\subset \BP \fg$ be the adjoint variety of $\fg$.
Let $X$ be a smooth projective variety with a family of minimal rational curves whose VMRT $\sC_x \in \BP T_x X$ at a general point $x\in X$ is isomorphic to $Y \subset \BP \fg$. Then the germ of  a minimal rational curve through $x$ is isomorphic to either the germ of a minimal rational curve in Example \ref{ex.flat} or
\begin{itemize} \item[(a)] a germ of a general line on a smooth hyperplane section of the Grassmannian ${\rm Gr}(3; \C^6) \subset \BP (\wedge^3 \C^6)$ for $\fg$  of type $\textsf{A}_2$; and
\item[(b)] a germ of a general minimal rational curve on the wonderful group compactification for $\fg$ of other types. \end{itemize} \end{theorem}

In (b), that  the VMRT-structures of wonderful group compactifications are $Y$-isotrivial is a result of Brion and Fu in \cite{BF}. The structure of hyperplane sections of  Grassmannians in (a) has been studied in arbitrary dimension by Bai, Fu and Manivel in \cite{BFM}. In fact, a result from \cite{BFM} shows that  when $Y \subset \BP \fg$ is the adjoint variety of $\fg$ of type $\textsf{A}_{\ell \geq 3}$,  there are examples of  $Y$-isotrivial cone structures with characteristic conic connections which are not locally symmetric. This shows   that the statement (b) in Theorem \ref{t.adjoint1} does not hold when $\fg$ is of type $\textsf{A}_{\ell \geq 3}$. We explain these in Example \ref{ex.A} below.

\medskip
Let us briefly comment on the proofs of the above theorems and the organization of the paper. The proofs of Theorems \ref{t.tensor1}, \ref{t.torsion1}, \ref{t.torsionfree1}   use the classical theory of principal connections on $G$-structures. We  review this classical theory in Section \ref{s.Gstructure}, including the local coordinates computation for the associated affine connections, which is required in the proof of Theorem \ref{t.torsion1}. In Section \ref{s.cone},  after reviewing the basic theory of cone structures and conic connections, we give the proofs of Theorems \ref{t.tensor1}, \ref{t.torsion1} and \ref{t.torsionfree1}.  In Section \ref{s.MS}, we review a result of \cite{MS} and explain how  it implies Theorem \ref{t.MS1},  when combined with Theorem \ref{t.torsionfree1}.  We review the geometric properties of the adjoint varieties in  Section \ref{s.adjoint} and use them to prove Theorem \ref{t.adjoint1} in Section \ref{s.Spencer}. The proof of Theorem \ref{t.Bianchi} given in Section \ref{s.Bianchi},  is purely algebraic, using the theory of roots and weights of  simple Lie algebras. Finally, we explain in Section \ref{s.vmrt} the application to the theory of minimal rational curves in  algebraic geometry, including the proof of Theorem \ref{t.vmrt1}.

\section{Principal connections on G-structures}\label{s.Gstructure}

In this section, we recall some standard results on G-structures and principal connections.
  See   Chapters II and III of \cite{KN},  Section 2 of \cite{Sc} and Chapter VII  of \cite{St} for more detailed presentation. Note that a principal connection is called just a connection (on a principal bundle) in many references. Since there are so many different types of connections used here, it seems better to call it a principal connection, following the convention in \cite{CS} and \cite{Sc}.

\begin{definition}\label{d.principal}
Fix a vector space $V$ of dimension $n$. Let $M$ be a  complex manifold of dimension $n$. \begin{itemize} \item[(i)] The {\em frame bundle} of $M$ is a principal ${\rm GL}(V)$-bundle  $\varpi: \Fr M  \rightarrow M$,  whose fiber at $x \in M,$ $$\Fr_x M := {\rm Isom}(V, T_x M)$$ is the set of linear isomorphisms from $V$ to $T_x M$. \item[(ii)] For a closed subgroup $G \subset {\rm GL}(V)$, a $G$-{\em structure} on $M$ (equivalently, a  G-structure on $M$ with {the structure group} $G$)  means a $G$-principal subbundle $\sP \subset \Fr M$.
\end{itemize}
Assume that a $G$-structure $\sP$ is given on $M$. Denote by $\fg \subset \fgl(V)= \End(V)$  the Lie algebra of $G \subset {\rm GL}(V)$, by $\pi: \sP \to M$  the restriction $\varpi|_{\sP}$ and by $T^{\pi} := {\rm Ker}({\rm d} \pi) \subset T \sP$  the vertical tangent bundle of $\pi$. \begin{itemize}
\item[(iii)]
For any point $\ep \in \sP$, we have a canonical identification of $\fg$ with
    $T^{\pi}_{\ep} \subset T_h \sP$ sending an element $\eta \in \fg $ to the vector $$\vec{\eta}_{\ep} := \frac{\rm d}{{\rm d} t}|_{t=0} (\ep \cdot \exp (t \eta)) \in T^{\pi}_{\ep},$$ where $\{\exp( t \eta) \in G \mid t \in \C\}$ is the 1-parameter subgroup of $G$ generated by $\eta$.  The resulting vector field $\vec{\eta}$ on $\sP$ is called the {\em fundamental vector field generated by} $\eta \in \fg.$
    \item[(iv)] The {\em soldering form} $\theta$ is a $V$-valued 1-form  on $\sP$  such that $\theta_{\ep}: T_{\ep} \sP \to V$ at the point $\ep \in \sP_x, x \in M$ is the composition
    $$T_{\ep} \sP \stackrel{\rd \pi}{\longrightarrow} T_x M \stackrel{\ep^{-1}}{\longrightarrow} V$$ where $\ep^{-1}$ is the inverse of the isomorphism $\ep\in \sP_x \subset {\rm Isom}(V, T_x M)$.
  \item[(v)] The fiber bundle $\sG := \sP \times_G \fg$   via the adjoint action of  $G$ on its Lie algebra $\fg$  is the {\em adjoint bundle} of the $G$-structure.   The tangent bundle $TM$ can be identified with $\sP \times_G V $ by the action of $G\subset {\rm GL}(V)$ on $V$. From $\fg \subset \End(V)$, the adjoint bundle $\sG$ can be naturally regarded as a vector subbundle of $\End(TM)$.
       \end{itemize}\end{definition}

The following is (2.5) in p. 310 of \cite{St}.

\begin{lemma}\label{l.St}
In Definition \ref{d.principal} (iii) and (iv),  for any $\ep \in \sP, w \in T_{\ep} \sP $ and $\eta \in \fg,$
$$ \rd \theta (\vec{\eta}, w) = - \eta \cdot \theta(w), $$ where $\eta \cdot $ means the action of $\eta \in \fg \subset \End (V)$ on $V$. \end{lemma}

We use the following definition of a principal connection and its torsion tensor, from p. 120 and p. 132 of \cite{KN}.

   \begin{definition}\label{d.connection}
   Let  $\sP \subset \Fr M$ be a $G$-structure as in Definition \ref{d.principal}.
   \begin{itemize} \item[(i)]
            A {\em principal connection} on $\sP$ is a vector subbundle $\sH \subset T \sP$  such that $T \sP = T^{\pi} \oplus \sH$ and ${\rm d} R_{g} (\sH )= \sH$ for all $g \in G$,  where ${\rm d} R_{g}: T\sP \to T\sP$  is the differential of the right action of $g$ on $\sP$. \end{itemize}
            Assume that a principal connection $\sH$ is given. \begin{itemize}
\item[(ii)]
            For the fiber $\sH_{\ep}$ at $\ep \in \sP_x$ and a vector $u \in T_x M$, denote by  $u^{\sH_{\ep}}$ the unique vector in $\sH_{\ep}$ such that ${\rm d} \pi (u^{\sH_{\ep}}) = u$.
            \item[(iii)] The {\em torsion tensor} of $\sH$ is a homomorphism  $\tau^{\sH} \in \Hom(\wedge^2 TM, TM)$  that satisfies
                $$\tau^{\sH}(u, w) = 2 \ \ep (\rd \theta (u^{\sH_{\ep}}, w^{\sH_{\ep}}))$$
                for any $\ep \in \sP_x$ and  any $u, w \in T_x M$.  \end{itemize}
\end{definition}

We skip the proof of the following easy lemma.

\begin{lemma}\label{l.iota}
Let $\sP$ be a $G$-structure on a complex manifold $M$ and let $\sG \subset \End(TM)$ be its adjoint bundle.  \begin{itemize} \item[(i)]
Each $\ep \in \sP_x \subset {\rm Isom}(V, T_x M)$ induces an isomorphism $\ep_*: \End(V) \stackrel{\cong}{\longrightarrow} \End(T_x M),$ the restriction of which gives  an isomorphism  $\iota_{\ep}: T^{\pi}_{\ep} \to \sG_x$  sending $ \vec{\eta}_{\ep} \in T^{\pi}_{\ep}$ generated by $\eta \in \fg \subset \End(V)$ to $$\iota_{\ep}(\vec{\eta}_{\ep}) :=  \ep_*(\eta) = \ep \circ \eta \circ \ep^{-1} \in \End(T_x M).$$
\item[(ii)]
 Denote by ${\rm d} R_g: T_{\ep} \sP \to T_{\ep \cdot g} \sP$  the differential of the right action by $g \in G$. Then $ \iota_{\ep} = \iota_{\ep \cdot g} \circ {\rm d}R_g.$ \end{itemize} \end{lemma}

\begin{proposition}\label{p.beta}
Let $\sP , \sG$ and $\iota_{\epsilon}, \epsilon \in \sP$ be as  in Lemma \ref{l.iota}. \begin{itemize} \item[(i)]
Let $\sH$ and $\widetilde{\sH}$ be two principal connections on $\sP$.
Note that for any  $\ep \in \sP_x$ and $u \in T_x M$, the difference $u^{\widetilde{\sH_{\ep}}} - u^{\sH_{\ep}} $ belongs to $T^{\pi}_{\ep}.$
Then there exists a homomorphism of vector bundles $\beta_{\sH,\widetilde{\sH}} \in \Hom(TM, \sG)$ such that for any  $\ep \in \sP_x$
and $u \in T_x M$, \begin{equation}\label{e.beta} \iota_{\ep}(u^{\widetilde{\sH}_{\ep}} - u^{\sH_{\ep}})  = \beta_{\sH,\widetilde{\sH}} (u)  \in \sG_x \subset \End(T_x M). \end{equation}  \item[(ii)] Conversely, for any homomorphism $\beta \in \Hom(TM, \sG)$ and a principal connection $\sH$ on $\sP$, there exists a unique principal connection $\widetilde{\sH}$ on $\sP$ such that $\beta = \beta_{\sH,\widetilde{\sH}}$. \end{itemize} \end{proposition}

\begin{proof}
To prove (i), just define $\beta_{\sH,\widetilde{\sH}}(u) \in \sG_x$ for each $u \in T_x M$ by  $$\beta_{\sH,\widetilde{\sH}}(u) := \iota_{\ep}(u^{\widetilde{\sH}_{\ep}} - u^{\sH_{\ep}}) \mbox{ for some } \ep \in \sP_x.$$
To check that this definition does not depend on the choice of $\ep \in \sP_x$,
it suffices to show that for any $g \in G$,
$$
\iota_{\ep}(u^{\widetilde{\sH}_{\ep}} - u^{\sH_{\ep}}) = \iota_{\ep \cdot g} (u^{\widetilde{\sH}_{\ep \cdot g}} - u^{\sH_{\ep\cdot g}}).$$
Since $$u^{\widetilde{\sH}_{\ep \cdot g}} - u^{\sH_{\ep \cdot g}} = {\rm d} R_g (u^{\widetilde{\sH}_{\ep}} - u^{\sH_{\ep}}),$$ this  follows from Lemma \ref{l.iota}. Thus we obtain $\beta_{\sH,\widetilde{\sH}} \in \Hom(TM,\sG)$  satisfying (\ref{e.beta}).

To prove (ii), we define a vector subbundle $\widetilde{\sH} \subset T \sP$ complementary to $T^{\pi}$ as follows. Note that whenever we have a subspace $H \subset T_{\ep} \sP$ at $\ep \in \sP_x$ complementary to $T^{\pi}_{\ep}$, we have an isomorphism $T_x M \to H$
sending $u \in T_x M$ to a vector $u^H \in H$.  Given a principal connection $\sH$ and $\beta \in \Hom(TM, \sG)$, we can find  a unique subspace $\widetilde{\sH}_{\ep} \subset T_{\ep} \sP$ complementary to $T^{\pi}_{\ep}$  at each $\ep \in \sH$ such that $$u^{\widetilde{\sH}_{\ep}} = u^{\sH_{\ep}} + \iota_{\ep}^{-1}(\beta(u)) \mbox{ for any } u \in T_x M.$$
This determines a subbundle $\widetilde{\sH} \subset T \sP$ complementary to $T^{\pi}.$
It is easy to see $\rd R_g (\widetilde{\sH}) = \widetilde{\sH}$ from Lemma \ref{l.iota}, which proves that $\widetilde{\sH}$ is a principal connection satisfying (\ref{e.beta}).
\end{proof}

 \begin{definition}\label{d.Spencer}
 Let $\sP$ be a $G$-structure on a complex manifold $M$. \begin{itemize}
\item[(i)] For the Lie algebra $\fg \subset \End(V)$ of the subgroup $G \subset {\rm GL}(V)$, the {\em Spencer homomorphism}
$\partial: \Hom(V, \fg) \to \Hom(\wedge^2 V, V)$ is defined
by $$\partial h (v_1, v_2) := h(v_1) \cdot v_2 - h(v_2) \cdot v_1 $$ for any $v_1, v_2 \in V$ and $h \in \Hom(V, \fg),$ where $h(v_i) \cdot v_j$ means the action of $h(v_i) \in \fg \subset \End(V)$ on $v_j$. \item[(ii)]
The Spencer homomorphism induces  a homomorphism $$\partial: \Hom(TM, \sG) \to \Hom(\wedge^2 TM, TM)$$
defined by
$$\partial \beta(u, w) = \beta(u) \cdot w - \beta(w) \cdot u, $$ for any $u, w \in T_x M$ and $\beta \in \Hom(TM, \sG)$, where $\beta(u) \cdot w$ (resp. $\beta(w) \cdot u$) means the action of $\beta(u)$ (resp. $ \beta(w)$) in  $\sG_x \subset \End(T_x M)$ on $w$ (resp. $u$).
\end{itemize}
\end{definition}

\begin{lemma}\label{l.Stein}
In Definition \ref{d.Spencer}, let ${\rm Im}(\p) \subset \Hom(\wedge^2 V, V)$ be the image of the Spencer homomorphism and let $\sP \times_G {\rm Im}(\p) \subset Hom(\wedge^2 TM, TM)$ be the corresponding vector subbundle. Then for any point $o \in M$, there exists a neighborhood $O$ such that $$\p (\Hom(TO, \sG|_O)) = H^0(O, \sP|_O \times_G {\rm Im}(\p)).$$\end{lemma}

\begin{proof} Consider the exact sequence of vector bundles on $M$:
$$ 0 \longrightarrow \sP \times_G {\rm Ker}(\p) \longrightarrow \sP \times_G \Hom(V, \fg)
\longrightarrow \sP \times_G {\rm Im}(\p) \longrightarrow 0.$$
Choose a Stein neighborhood $O$ of $o \in M$
such that $H^1(O, \sP \times_G {\rm Ker}(\p)) =0.$
Then $$ \Hom(TO, \sG|_O)  = H^0(O, \sP|_O \times_G \Hom(V, \fg)) \stackrel{\p}{\longrightarrow} H^0(O, \sP|_{O} \times_G {\rm Im}(\p))$$ is surjective.
\end{proof}

\begin{proposition}\label{p.torsion}
Let $\sH, \widetilde{\sH}$ and $\beta_{\sH, \widetilde{\sH}}$ be as in Proposition \ref{p.beta} (i).

Then $\tau^{\sH} - \tau^{\widetilde{\sH}} = 2 \ \partial \beta_{\sH, \widetilde{\sH}}.$ \end{proposition}

\begin{proof} For $u, w \in T_x M$ and $\ep \in \sP_x$,
\begin{eqnarray*} \lefteqn{ \frac{1}{2}\ep^{-1}(\tau^{\sH}(u,w) - \tau^{\widetilde{\sH}}(u,w) )} \\   && =
\rd \theta (u^{\sH_{\ep}}, w^{\sH_{\ep}})- \rd \theta(u^{\widetilde{\sH}_{\ep}}, w^{\widetilde{\sH}_{\ep}}) \\ && =
\rd \theta ((u^{\sH_{\ep}} - u^{\widetilde{\sH}_{\ep}}), w^{\sH_{\ep}}) + \rd \theta ( u^{\widetilde{\sH}_{\ep}}, ( w^{\sH_{\ep}} - w^{\widetilde{\sH}_{\ep}})) \\ && = \rd \theta(- \iota_{\ep}^{-1}(\beta_{\sH, \widetilde{\sH}}(u)), w^{\sH_{\ep}}) + \rd \theta(u^{\widetilde{\sH}_{\ep}}, - \iota_{\ep}^{-1}(\beta_{\sH, \widetilde{\sH}}(w))) \\ && = \ep_*^{-1}(\beta_{\sH, \widetilde{\sH}}(u)) \cdot \theta(w^{\sH_{\ep}}) - \ep_*^{-1}(\beta_{\sH, \widetilde{\sH}}(w)) \cdot \theta(u^{\widetilde{\sH}_{\ep}}),
 \end{eqnarray*} where the last equality uses Lemma \ref{l.St}.
 Thus $\tau^{\sH}(u,w) - \tau^{\widetilde{\sH}}(u,w) $ equals to
\begin{eqnarray*} & & 2 \, \ep (\ep_*^{-1}(\beta_{\sH, \widetilde{\sH}}(u)) \cdot \ep^{-1}(w) - \ep_*^{-1}(\beta_{\sH, \widetilde{\sH}}(w) \cdot \ep^{-1}(u))) \\ &=& 2 \, ( \beta_{\sH, \widetilde{\sH}}(u) \cdot w - \beta_{\sH, \widetilde{\sH}}(w) \cdot u) \\ &=& 2 \, \p \beta_{\sH, \widetilde{\sH}} (u, w). \end{eqnarray*} \end{proof}

\begin{definition}\label{d.affine}
Let $\sH$ be a principal connection on a $G$-structure $\sP \subset \Fr M.$
\begin{itemize} \item[(i)] The  quotient morphism  $\sP \times V \to \sP \times_G V = TM$ by the $G$-action $(\ep, v) \cong (\ep \cdot g^{-1}, g \cdot v)$  sends
$\sH \subset T\sP$ to a vector subbundle $\tH \subset T(TM),$ called the {\em affine connection} on $M$ induced by $\sH$,  such that $T(TM) = T^{\psi} \oplus  \tH$ where $T^{\psi}$  denotes the vertical tangent bundle of the natural projection $\psi: TM \to M$. \item[(ii)] For the fiber $\tH_w$ at $w \in T_x M$ and a vector $u \in T_x M$, denote by  $u^{\tH_w}$ the unique vector in $\tH_w$ such that  $\rd \psi(u^{\tH_w}) = u$.
\item[(iii)] Let $a: \Delta \to M$ be any arc in $M$ with $a(0) = x$. For a vector $u \in T_xM,$ we have a unique  arc $a^{\sharp}: \Delta \to TM$ such that $$a^{\sharp}(0) = u, \ \psi \circ a^{\sharp} = a, \ \rd a^{\sharp} (T \Delta) \subset \tH.$$ This arc $a^{\sharp}$ or its image $a^{\sharp}(\Delta) \subset TM$ is called the $\sH$-{\em parallel transport} of $u$ along the arc $a$.
\item[(iv)] The {\em geodesic flow} of the principal connection $\sH$ is the vector field  $\vec{\gamma}^{\sH}$ on $TM$ defined by $\vec{\gamma}^{\sH}_w := w^{\tH_w}$ for each $w \in TM$.
    \item[(v)] The $\psi$-image of an integral curve of $\vec{\gamma}^{\sH}$ is called an {\em $\sH$-geodesic}. \end{itemize} \end{definition}

The following lemma is immediate.

\begin{lemma}\label{l.descend}
In Definition \ref{d.affine}, let $A \subset V$ be a submanifold preserved by the $G$-action on $V$ and let $\sA \subset TM$ be the submanifold given by $$ \sA:=  \sP \times_G A \subset \sP \times_G V = TM.$$ Then $\tH_a \subset T_a \sA$ for any  $a \in \sA$, namely, the affine connection $\tH$  induced by $\sH$ is tangent to $\sA$.
 \end{lemma}

 The following is  a direct consequence of Proposition \ref{p.beta} (i).

 \begin{proposition}\label{p.beta2}
 Let $\sP \subset \Fr M$ be a $G$-structure on a complex manifold $M$ and let
$\sH$, $\widetilde{\sH}$ and $\beta_{\sH, \widetilde{\sH}} \in \Hom( TM, \sG)$ be as in Proposition \ref{p.beta} (i). Let $\tH$ and $\widetilde{\tH}$ be the  affine connections on $M$ induced by $\sH$ and $\widetilde{\sH}$, respectively. Then$$ u^{\widetilde{\tH}_w} - u^{\tH_w} =  \beta_{\sH, \widetilde{\sH}}(u) \cdot w \in T_x M$$ for any $x \in M$ and $u, w \in T_x M$.
\end{proposition}

 \begin{definition}\label{d.nabla}
Let $\sH \subset T \sP$ be a principal connection on a $G$-structure and let $\tH \subset T(TM)$ be the induced affine connection.
\begin{itemize} \item[(i)] For each $w \in T_x M$, we have the projection  $T_{w} (TM)  \to T^{\psi}_{w}$  with the kernel $\tH_{w}$ and  the canonical identification $T^{\psi}_{w} = T_x M.$ Their composition defines a surjective homomorphism $p^{\sH}_{w} : T_{w} (TM) \to T_x M$.  This determines a morphism $p^{\sH}:  T(TM) \to TM$ with ${\rm Ker}(p^{\sH}) = \tH$.  \item[(ii)]
For any vector fields $\vec{u}, \vec{w}$  on an open subset $U \subset M$, define the vector field $\nabla^{\sH}_{\vec{u}} \vec{w}$ on $U$ as follows.  Let $\sigma: U \to TU$ be the section of $\psi: TU \to U$ given by $\vec{w}$ and define $ \nabla^{\sH}_{\vec{u}} \vec{w}$   as the section of $TU$ given by $p^{\sH} ({\rm d} \sigma(\vec{u})).$ The operator $\nabla^{\sH}$ which sends a pair of vector fields $\vec{u}, \vec{w}$ on any open subset $U \subset M$ to a vector field $\nabla^{\sH}_{\vec{u}} \vec{w}$ on $U$  is the {\em covariant derivative} induced by $\sH$.
 \end{itemize} \end{definition}

The following is standard (see  Theorem 5.1 in Chapter III of \cite{KN}).

\begin{proposition}\label{p.torsionnabla}
In Definition \ref{d.affine}, the torsion tensor $\tau^{\sH}$ satisfies
 $$\tau^{\sH}(v,w ) = \nabla^{\sH}_{\vec{v}} \vec{w} - \nabla^{\sH}_{\vec{w}} \vec{v} - [\vec{v}, \vec{w}] $$ for any $v, w \in T_x M$ and any local vector fields $\vec{v}, \vec{w}$ with $\vec{v}_x = v$ and $\vec{w}_x = w.$
 \end{proposition}

 \begin{definition}\label{d.curvature}
 Let $\nabla^{\sH}$ be as in Definition \ref{d.nabla}.
\begin{itemize} \item[(i)] The {\em curvature} of the principal connection $\sH$ is a homomorphism
$R^{\sH} \in \Hom(\wedge^2 TM, \End(TM))$ defined by
$$R^{\sH}(v,w) \cdot u = \nabla^{\sH}_{\vec{v}} (\nabla^{\sH}_{\vec{w}} \vec{u}) - \nabla^{\sH}_{\vec{w}}(\nabla^{\sH}_{\vec{v}}\vec{u}) - \nabla^{\sH}_{[\vec{v}, \vec{w}]} \vec{u}$$ for any $v, w, u \in T_x M$ and any local vector fields $\vec{v}, \vec{w}, \vec{u}$ with $\vec{v}_x = v, \vec{w}_x = w$ and $\vec{u}_x = u$.
\item[(ii)] The principal connection $\sH$ is {\em torsion-free} if $\tau^{\sH} =0$. The $G$-structure $\sP$ is {\em torsion-free} if  each point $x \in M$ admits a neighborhood $x \in U \subset M$ with a torsion-free principal connection on $\sP|_U$.
\item[(iii)]
The principal connection $\sH$ is said to be {\em locally flat}  if $\tau^{\sH} = R^{\sH} =0.$ The $G$-structure $\sP$ is {\em locally flat} if  each point $x \in M$ admits a neighborhood $x \in U \subset M$ with a locally flat principal connection on $\sP|_U$.
    \item[(iv)] Note that the covariant derivative $\nabla^{\sH}_{\vec{w}}$ can be extended to operate on any tensor fields. The principal connection $\sH$ is said to be {\em locally symmetric} if  $\tau^{\sH} =0$ and $\nabla_{\vec{w}}^{\sH} R^{\sH} =0$ for any local vector field $\vec{w}$. The $G$-structure $\sP$ is {\em locally symmetric} if  each point $x \in M$ admits a neighborhood $x \in U \subset M$ with a locally symmetric principal connection on $\sP|_U$.
\end{itemize} \end{definition}

\begin{definition}\label{d.Bianchi}
Let $\fg \subset \End(V)$ be a Lie subalgebra. For each $h \in \Hom(\wedge^2 V, \fg)$, define $h^{\sharp} \in \Hom(\wedge^3 V, V)$ by $$h^{\sharp}(u,v,w) = h(u,v) \cdot w + h(v,w) \cdot u + h(w,u) \cdot v$$ for all $u, v, w \in V,$ where $h(u,v) \cdot w$ denotes the action of $h(u,v) \in \End(V)$ on $w \in W$. Define
$$\BK (\fg ) := \{ h \in \Hom(\wedge^2 V, \fg) \mid h^{\sharp} =0\}.$$ This subspace $\BK(\fg)$ of $\Hom(\wedge^2 V, \fg)$ is called the {\em space of formal curvature maps } of $\fg \subset \End(V)$.
\end{definition}

The following is standard (see  Theorem 5.3 in Chapter III of \cite{KN} or  page 16 of \cite{Sc}).

\begin{proposition}\label{p.BK}
Let $\sH$ be a torsion-free  principal connection on a G-structure $\sP$ on a complex manifold $M$ and let $\fg \subset \End(V)$ be the Lie algebra of the structure group of $\sP$. Then the curvature $R^{\sH} \in \Hom(\wedge^2 TM, \End(TM))$ belongs to $\Hom(\wedge^2 TM, \sG)$ and satisfies (called the first Bianchi identity)
$$R^{\sH}(u,v) \cdot w + R^{\sH}(v,w) \cdot u + R^{\sH}(w,u) \cdot v =0.$$ In other words, it is a holomorphic section of the vector subbundle
$$\sP \times_G \BK (\fg) \ \subset \sP \times_G \Hom(\wedge^2 V, \fg) = \wedge^2 T^*M \otimes \sG.$$ \end{proposition}

We have the following coordinate expressions from   Chapter III, Proposition 7.6, Proposition 7.8 of \cite{KN} and the equation (4.13) of \cite{Lee}.

\begin{proposition}\label{p.coordi}
In Definition \ref{d.nabla},
suppose that we have a coordinate system $(x^1, \ldots, x^n), n = \dim M,$  on an open subset $U \subset M$.  Define holomorphic functions $\Gamma^{k}_{ij}$ on $U$ for $1 \leq i, j, k \leq n$,   by $$\nabla^{\sH}_{\frac{\p}{\p x^i}} \frac{\p}{\p x^j} =  \sum_{k =1}^n \Gamma^k_{ij} \frac{\p}{\p x^k} \mbox{ for all } 1 \leq i, j \leq n.$$
\begin{itemize} \item[(i)] The torsion $\tau^{\sH}$ satisfies
$$ \tau^{\sH}(\frac{\p}{\p x^i}, \frac{\p}{\p x^j}) = \sum_{k=1}^n (\Gamma^k_{ij} - \Gamma^k_{ji}) \frac{\p}{\p x^k} $$  for all $ 1 \leq i, j \leq n.$  \item[(ii)]
Choose a coordinate system on $TM$ by $z^i = \psi^* x^i$ and $y^i = {\rm d} x^i$ for $1\leq i \leq n$. Then the geodesic flow of $\sH$ is  given by $$\vec{\gamma}^{\sH} = \sum_{k =1}^n y^k \frac{\p}{\p z^k} - \sum_{i,j,k =1}^n \Gamma^k_{ij} y^i y^j \frac{\p}{\p y^k}. $$
\item[(iii)]
Let $\{a(t)= (a^1(t), \ldots, a^n(t)) \mid t \in \Delta\}$ be an arc in $U$ and let $$\{ a^{\sharp}(t) = \sum_{j=1}^n c^j(t) \frac{\p}{\p x^j}\in T_{a(t)} M \}$$ be the $\sH$-parallel transport of $a^{\sharp}(0) = \sum_{j=1}^n c^j(0) \frac{\p}{\p x^j}$ along the arc $a$. Then
$$\frac{\rd}{\rd t} c^k(t) = -  \sum_{i, j =1}^n\Gamma^k_{i,j =1}(a(t))  \frac{\rd a^i(t)}{\rd t}  c^j(t)  $$ for all $t \in \Delta$ and $1 \leq k \leq n$.  \end{itemize} \end{proposition}

\section{Conic connections on cone structures}\label{s.cone}
\begin{definition}\label{d.cone} Let $M$ be a complex manifold and let $\BP TM$ be the projectivization of its tangent bundle.
\begin{itemize}
\item[(i)] A {\em cone structure} on $M$ is a closed submanifold $\sC \subset \BP TM$ which is submersive over $M$. Denote by $\phi: \sC \to M$ the projection to $M$. Then the fiber $\sC_x := \phi^{-1}(x) \subset \BP T_x M$ is a projective submanifold  at every point $x \in M.$ Write $L$ for the tautological line bundle on $\sC$ whose fiber $L_{\alpha}$ at $\alpha \in \sC_x$ is $\widehat{\alpha} \subset T_x M$.
    \item[(ii)] For each $\alpha \in \sC$ in (i), let $\sJ_{\alpha} \subset T_{\alpha} \sC$ be the subspace $({\rm d} \phi)^{-1}(\widehat{\alpha})$. The vector subbundle $\sJ \subset T\sC$ defined by $\{\sJ_{\alpha} \mid \alpha \in \sC\}$ is the {\em tautological distribution} on $\sC$.
        \item[(iii)] In (ii), we have an exact sequence $$0 \to T^{\phi} \to \sJ \to L \to 0$$ on $\sC$. A line subbundle $\sF \subset \sJ$ is a {\em conic connection} of the cone structure $\sC$ if it splits this exact sequence, namely, there is a direct sum decomposition $\sJ = T^{\phi} \oplus \sF$ inducing an isomorphism of line bundles $\sF \cong L$.
            \item[(iv)] In (ii), let $\sE_{\alpha} \subset T_{\alpha} \sC$ be the subspace $({\rm d} \phi)^{-1}(T_{\alpha} \widehat{\sC}_x)$. The collection $\{\sE_{\alpha} \mid \alpha \in \sC \}$ determines a vector subbundle $\sE \subset T\sC$.  A conic connection $\sF \subset \sJ$ is a {\em characteristic conic connection}  if the sheaves $\sO(T^{\phi}), \sO(\sF)$ and $\sO(\sE)$ of vector fields on $\sC$  given by local sections of $T^{\phi}, \sF$ and $\sE,$ respectively,  satisfy $$[[ \sO(T^{\phi}), \sO(\sF)], \sO(\sF)] \subset \sO(\sE).$$
            \end{itemize} \end{definition}

\begin{definition}\label{d.isotrivial}
In Definition \ref{d.cone}, fix a vector space $V$ with $\dim V = \dim M$ and a  projective submanifold $Z \subset \BP V$.
\begin{itemize} \item[(i)] A cone structure $\sC \subset \BP TM$ is {\em $Z$-isotrivial} if for every $x \in M$, there is a linear isomorphism $\ep \in {\rm Isom}(V, T_x M)$ such that $\ep(\widehat{Z}) = \widehat{\sC}_x$.  A cone structure $\sC$ is {\em isotrivial} if it is $Z$-isotrivial for some $Z \subset \BP V$.
\item[(ii)] Given a $Z$-isotrivial cone structure $\sC \subset \BP TM$, define $\sP_x \subset \Fr_x M$ by
$$\sP_x := \{ \ep \in {\rm Isom}(V, T_x M) \mid \ep(\widehat{Z}) = \widehat{\sC}_x\}.$$
Then $\sP \subset \Fr M$ is a G-structure on $M$ with the group $${\rm Aut}(\widehat{Z}) := \{ g \in {\rm GL}(V) \mid g \cdot Z = Z \}.$$ This is called the {\em associated} G-{\em structure} of the isotrivial cone structure $\sC \subset \BP TM$.
\end{itemize} \end{definition}

We omit the easy proof of the following.

\begin{proposition}\label{p.isotrivial}
Let $Z \subset \BP V$ be a projective submanifold.  Let $\sP \subset \Fr M$ be a G-structure on a complex manifold $M$ with the structure group ${\rm Aut}(\widehat{Z}) \subset {\rm GL}(V).$
Then the images of $Z$ under elements of $$\{ \sP_x \subset {\rm Isom}(V, T_x M) \mid x \in M\}$$ determine a $Z$-isotrivial cone structure $\sC \subset \BP TM$ such that $\sP$ is the associated G-structure of $\sC \subset \BP TM$ in the sense of Definition \ref{d.isotrivial}. \end{proposition}

\begin{proposition}\label{p.inducedconicconnection}
Let $\sC \subset \BP TM$ be a $Z$-isotrivial cone structure and let $\sP \subset \Fr M$ be the associated G-structure. Let $\sH \subset T \sP$ be a principal connection.
\begin{itemize}
\item[(i)] For any arc $a:\Delta \to M$ and any point $c \in \check{\sC} \subset TM$, the $\sH$-parallel transport of $c$ along $a$ is contained in $\check{\sC}$.
 \item[(ii)] The geodesic flow $\vec{\gamma}^{\sH} \subset T(TM)$ is tangent to $\check{\sC} \subset TM \setminus O_M$ and the $\C^{\times}$-projection $ \check{\sC} \to \sC \subset \BP TM$ sends $\vec{\gamma}^{\sH}|_{\check{\sC}}$  to a conic connection on $\sC$. \end{itemize} \end{proposition}

\begin{proof}
The affine connection $\tH \subset T(TM)$ induced by $\sH$ is tangent to $\check{\sC}$ by Lemma \ref{l.descend}. Both (i) and (ii) follow from this. \end{proof}

\begin{definition}\label{d.inducedconic}
The conic connection on $\sC \subset \BP TM$  in Proposition  \ref{p.inducedconicconnection} given by the image of $\vec{\gamma}^{\sH}|_{\check{\sC}}$ is the {\em conic connection induced by} $\sH$ and denoted by $\sF^{\sH}.$  \end{definition}

\begin{theorem}\label{t.CtoP}
Let $Z \subset \BP V$ be a closed submanifold. Let  $\aut(\widehat{Z}) \subset \End(V)$ be the Lie algebra of $\Aut(\widehat{Z}) \subset {\rm GL}(V)$. Suppose that
\begin{itemize} \item[(i)] $H^0(Z, \sO(1)) = V^*$; and \item[(ii)] the homomorphism $$\aut(\widehat{Z}) \otimes V^* \to H^0(Z, TZ \otimes \sO(1))$$ induced by the restriction homomorphism $\aut(\widehat{Z}) \to H^0(Z, TZ)$ is surjective. \end{itemize}
Then any conic connection on a $Z$-isotrivial  cone structure is locally a conic connection induced by a principal connection on the associated G-structure. More precisely,  for any $Z$-isotrivial cone structure  $\sC \subset \BP TM$ on a complex manifold $M$ equipped with a conic connection  $\sF \subset \sJ \subset T \sC$,  each point $o \in M$ admits a neighborhood $o \in O \subset M$ such that the restriction $\sP|_O$ of the associated G-structure $\sP$  to $O$ has a principal connection $\sH$  satisfying $\sF|_{\phi^{-1}(O)} = \sF^{\sH}$. \end{theorem}

    \begin{proof}
  We can choose a neighborhood $o \in O \subset M$  such that  $\sP|_O$ admits   a principal connection $\widetilde{\sH} \subset T \sP|_{O}.$  Then we have the induced conic connection $\sF^{\widetilde{\sH}}$ on $\sC|_O$. Since both $\sF|_O$ and $\sF^{\widetilde{\sH}}$ are splittings of the exact sequence in Definition \ref{d.cone} (iii), their difference is an element $f \in H^0(\sC|_O, T^{\phi} \otimes L^{-1}).$ For each point $x \in O$, the restriction
$$  f|_{\sC_x} \in H^0(\sC_x, T \sC_x \otimes \sO(1)) \cong H^0(Z, TZ \otimes \sO(1))$$ comes from an element of $$\Hom(T_x M, \sG_x) = \aut(\widehat{\sC}_x) \otimes T^*_x M \cong \aut(\widehat{Z}) \otimes H^0(Z, \sO(1)),$$ by the assumptions (i) and (ii).  It follows that, after replacing $O$ by a smaller neighborhood if necessary,  we can find a homomorphism $\beta \in \Hom (TO, \sG|_O )$ that is sent to $f$ by the natural homomorphism
$$ \Hom(TO, \sG|_O) \to      H^0(\sC|_O, T^{\phi} \otimes L^{-1}).$$
 By Proposition \ref{p.beta} (ii), we obtain a principal connection $\sH$ on $\sP|_O$ satisfying $\sF = \sF^{\sH}$. \end{proof}

\begin{definition}\label{d.Xi}
For a projective submanifold $Z \subset \BP V,$ define the subspace $\Xi_Z \subset \Hom(\wedge^2 V, V)$ by
$$\Xi_Z := \{ \sigma \in \Hom(\wedge^2 V, V) \mid \sigma(v, T_v \widehat{Z}) \in T_v \widehat{Z} \mbox{ for any } v \in \widehat{Z} \}.$$ \end{definition}

\begin{theorem}\label{t.Xi}
Fix a projective submanifold $Z \subset \BP V$. Let $\sC \subset \BP TM$ be a $Z$-isotrivial cone structure on a complex manifold $M$ and let $\sP \subset \Fr M$ be the associated G-structure. Let $\sH \subset T \sP$ be a principal connection. We have  the induced conic connection $\sF^{\sH} \subset T \sC$ and the torsion $$\tau^{\sH} \in \Hom(\wedge^2 TM, TM) = H^0(M, \sP \times_G \Hom( \wedge^2 V, V)).$$ Then $\sF^{\sH}$ is a characteristic conic connection if and only if $\tau^{\sH}$ has values in $\Xi_Z$, namely, it belongs to the subspace
$$H^0(M, \sP \times_G \Xi_Z) \subset H^0(M, \sP \times_G \Hom(\wedge^2 V, V)).$$ \end{theorem}

\begin{proof}
Since the problem is local, we may replace $M$ by an open neighborhood of a base point $o \in M$,  equipped with a coordinate system $(x^1, \ldots, x^n)$. We use the notation from Proposition \ref{p.coordi}. Let us simplify the notation by using Einstein summation convention and writing $\nabla$ for the covariant derivative $\nabla^{\sH}$ induced by the principal connection $\sH$.

Choose a straight arc  through $o \in M$ $$a(t)  = (a^1 t, \ldots, a^n t) \in M,  \ t \in \Delta$$ in the direction of $(a^1, \ldots, a^n) \in \C^n$.
 Fix a point $b_o \in \check{\sC}_o$ and let $$\{ b(t) = c^j(t) \frac{\p}{\p x^j} \in T_{a(t)} M \mid t \in \Delta, b(0) = b_o\}$$ be the $\sH$-parallel transport of $b_o$ along the arc $a(t)$. By Proposition \ref{p.coordi} (iii), we have  \begin{equation}\label{e.parallel}
\frac{{\rm d}}{{\rm d} t} c^k(t) = -  \Gamma^k_{ij}  a^i c^j(t)  \mbox{ for all } t \in \Delta \mbox{ and } 1 \leq k \leq n. \end{equation}
  In the coordinate system $(z^i = \psi^* x^i, y^i = \rd x^i), 1 \leq i \leq n$ on $TM$, as in Proposition \ref{p.coordi},  the arc $b(t)$ is given by $$\{ z^k(t) = a^k t, \ y^k = c^k(t), \ t \in \Delta\}.$$
  Since the $\sH$-parallel transport must be contained in $\check{\sC}$, the tangent vector of $b(t)$ at $t=0$ satisfies
  $$\frac{\rd b(t)}{\rd t}|_{t=0} = a^i \frac{\p}{\p z^i} + \frac{\rd c^k}{\rd t}(0) \frac{\p}{\p y^k} \ \in T_{b_o} \check{\sC}.$$
Thus (\ref{e.parallel}) implies that \begin{equation}\label{e.tangent}
 a^i \frac{\p}{\p z^i} - \Gamma^k_{ij}(o) a^i c^j(0) \frac{\p}{\p y^k} \ \in T_{b_o} \check{\sC}.  \end{equation}

 Let $\vec{v} =  v^{\ell} \frac{\p}{\p y^{\ell}}$ be a section of $T^{\phi}$ in a neighborhood of $b_o$ in $TM$ which is tangent to $\check{\sC}$
and let $$\vec{f} =  y^k \frac{\p}{\p z^k} -  \Gamma^k_{ij} y^i y^j \frac{\p}{\p y^k}$$ be the geodesic flow of the principal connection $\sH$ from Proposition \ref{p.coordi} (ii). By a direct computation, \begin{equation}\label{e.vf}
[\vec{v}, \vec{f}] = v^i \frac{\p}{\p z^i} - \Gamma^k_{ij} v^i y^j \frac{\p}{\p y^k} + (\Gamma^{\ell}_{ij} y^i y^j \frac{\p v^k}{\p y^{\ell}} - \Gamma^k_{i \ell} y^i v^{\ell} - y^{\ell} \frac{\p v^k}{\p z^{\ell}}) \frac{\p}{\p y^k}. \end{equation}
From (\ref{e.tangent}), the first two terms on the right hand side  is tangent to $\check{\sC}$.
Thus we have $$ \clubsuit :=
(\Gamma^{\ell}_{ij} y^i y^j \frac{\p v^k}{\p y^{\ell}} - \Gamma^k_{i \ell} y^i v^{\ell} - y^{\ell} \frac{\p v^k}{\p z^{\ell}}) \frac{\p}{\p y^k}|_{b_o}  \ \in \  T_{b_o} \check{\sC}_o. $$
A direct computation using (\ref{e.vf}) shows that $[[\vec{v}, \vec{f}], \vec{f}],$ modulo local sections of $T^{\phi},$ equals to
\begin{equation}\label{e.modTphi}
( -2 y^{\ell} \frac{\p v^k}{\p z^{\ell}} + 2 \Gamma^{\ell}_{ij} y^i y^j \frac{\p v^k}{\p y^{\ell}} - \Gamma^k_{\ell j} v^{\ell} y^j - \Gamma^k_{\ell j} y^{\ell} v^j) \frac{\p}{\p z^k}. \end{equation}
By Definition \ref{d.cone} (iv), the vector field $\vec{f}$ is sent to a local section of a characteristic conic connection if and only if $[[\vec{v}, \vec{f}], \vec{f}]$ is sent to a local section of $\sE$ for all local section $\vec{v}$ of $T^{\phi}$ tangent to $\check{\sC}$. Thus $\sF^{\sH}$ is a characteristic conic connection if and only if for any choice of local section $\vec{v}$ of $T^{\phi}$ tangent to $\check{\sC}$,  (\ref{e.modTphi}) at  $b_o \in \check{\sC}$ is sent by ${\rm d} \phi$ to a vector in $T_o M$ tangent to $\check{\sC}_o \subset T_o M$.
This is equivalent to saying that $$ \spadesuit := ( -2 y^{\ell} \frac{\p v^k}{\p z^{\ell}} + 2 \Gamma^{\ell}_{ij} y^i y^j \frac{\p v^k}{\p y^{\ell}} - \Gamma^k_{\ell j} v^{\ell} y^j - \Gamma^k_{\ell j} y^{\ell} v^j) \frac{\p}{\p y^k}|_{b_o} $$ is tangent to $\check{\sC}_o$ at $b_o$. Note \begin{eqnarray*}  -2 \clubsuit + \spadesuit &=& (\Gamma^k_{\ell j} y^{\ell} v^j - \Gamma^k_{\ell j} v^{\ell} y^j) \frac{\p}{\p y^k}|_{b_o}\\ & = & (\Gamma^k_{\ell j} - \Gamma^k_{j \ell}) c^{\ell}(0) v^j(b_o) \frac{\p}{\p y^k}.\end{eqnarray*}
Thus $\sF^{\sH}$ is a characteristic connection if and only if $(\Gamma^k_{\ell j} - \Gamma^k_{j \ell})  c^{\ell} (0) v^j\frac{\p}{\p y^k}$ is tangent to $\check{\sC}_o$ for any $b_o \in \check{\sC}_o$ and $v^k \frac{\p}{\p y^k}|_{b_o} \in T_{b_o} \check{\sC}_o$. By Proposition \ref{p.coordi} (i), this is equivalent to saying that $\tau^{\sH}$ takes values in $\Xi_Z$.
\end{proof}

\begin{theorem}\label{t.Spencer}
Let $Z \subset \BP V$ be a projective submanifold satisfying the following three conditions:
\begin{itemize} \item[(i)] $H^0(Z, \sO(1)) = V^*$; \item[(ii)] the homomorphism $\aut(\widehat{Z}) \otimes H^0(Z, \sO(1)) \to H^0(Z, TZ(1))$ is surjective; and \item[(iii)] $\Xi_Z \subset \partial( \Hom(V, \aut(\widehat{Z})))$ in the notation of Definitions \ref{d.Spencer} and \ref{d.Xi}.  \end{itemize}   If a $Z$-isotrivial cone structure $\sC \subset \BP TM$ on a complex manifold $M$ admits  a characteristic conic connection, then the associated G-structure $\sP \subset \Fr M$ is torsion-free. \end{theorem}

\begin{proof}
Since the statement is local, we may  replace $M$ by a suitable open subset whenever it is necessary.
The conditions (i) and (ii) imply by Theorem \ref{t.CtoP}, after replacing $M$ by a suitable open subset, that there exists a principal connection $\sH$ on $\sP$ such that the induced  conic connection $\sF^{\sH}$ on $\sC$ is a characteristic conic connection. By  Theorem \ref{t.Xi}, the torsion tensor $\tau^{\sH}$ belongs to $H^0(M, \sP \times_G \Xi_Z)$. By the condition (iii),  $$H^0(M, \sP \times_G \Xi_Z) \subset H^0(M, \sP \times_G {\rm Im}(\p)).$$ Replacing $M$ by a suitable open subset, we have $$ H^0(M, \sP \times_G {\rm Im}(\p))
= \p \Hom(TM, \sG)$$ from Lemma \ref{l.Stein}.  Thus $\tau^{\sH} = \p \beta$ for some $\beta \in \Hom(TM, \sG)$. Then by Proposition \ref{p.torsion}, we can find a principal connection $\widetilde{\sH}$ on $\sP$ with $\tau^{\widetilde{\sH}} =0$. \end{proof}

\section{Torsion-free  irreducible G-structures}\label{s.MS}

In this section, we look at the following special type of $Z \subset \BP V$.

\begin{definition}\label{d.sky}
A closed connected  submanifold  $Z \subset \BP V$ is a {\em highest weight variety}  if the group $G = \Aut(\widehat{Z}) \subset {\rm GL}(V)$ is a reductive group which acts irreducibly on $V$ and  $Z$ is the orbit of highest weight vectors, in other words, the unique closed orbit of the $G$-action on $\BP V$. \end{definition}

We recall the  following standard results. (i) is Bott's Theorem (Theorem IV of \cite{Bo}) and  (ii) is a direct consequence of (i).

\begin{lemma}\label{l.Bott}
Let $Z \subset \BP V$ be a highest weight variety.
\begin{itemize} \item[(i)]
Write $Z = G/P$ for the isotropy subgroup $P$ of a base point $o \in Z$. Regarding the quotient $G \to G/P$ as a $P$-principal bundle on $Z$, for each irreducible representation $E$ of $P,$ we have the associated vector bundle $G\times_P E$ on $Z$.
Then $H^0(Z, G \times_P E)$ is an irreducible representation of $G$.
\item[(ii)] The restriction homomorphism $\Sym^k V^* \to H^0(Z, \sO(k))$ is surjective for all $k \geq 0.$ \end{itemize} \end{lemma}

The following is a classical result of Cartan and Kobayashi-Nagano (see Table 5 of \cite{MS}).

\begin{proposition}\label{p.KN}
Let $Z \subset \BP V$ be a highest weight variety. Then the Spencer homomorphism $\p: \Hom(V, \aut(\widehat{Z})) \to \Hom (\wedge^2 V, V)$ is not injective if and only if $\Aut(\widehat{Z}) \subset {\rm GL}(V)$ is the isotropy representation of an irreducible Hermitian symmetric space of compact type. More explicitly, the Spencer homomorphism $\p$ is not injective if and only if $Z \subset \BP V$ is isomorphic to one of the following projective subvarieties. \begin{itemize}
\item[(I)] Segre variety $\BP^{a-1} \times \BP^{b-1} \subset \BP^{ab-1}$.
\item[(II)] The Grassmannian of $2$-dimensional subspaces in a vector space $\C^m$ of dimension $m$, under the Pl\"ucker embedding in $\BP (\wedge^2 \C^m)$.
\item[(III)] The second Veronese embedding of $\BP^{m-1}$ in $\BP (\Sym^2 \C^m)$.
\item[(IV)] The nonsingular quadric hypersurface $\Q^{n-2} \subset \BP^{n-1}$ for $n\geq 3$.
\item[(V)] The 10-dimensional highest weight variety  of the $16$-dimensional spin representation of $\fso_{10}.$
\item[(VI)] The 16-dimensional highest weight variety of the $27$-dimensional irreducible representation of the simple Lie algebra of type $\textsf{E}_6$.
\end{itemize}
\end{proposition}

\begin{definition}\label{d.Legendre}
A highest weight variety $Z \subset \BP V$ is called a {\em homogeneous Legendrian variety} (or a {\em subadjoint variety}) if there exists a symplectic form $\omega$ on $V$ with respect to which $\widehat{Z} \subset V$ is Lagrangian, namely,
\begin{itemize} \item[(1)] $\dim \widehat{Z} = \frac{1}{2} \cdot \dim V$ and
\item[(2)] $\omega(T_z \widehat{Z}, T_z \widehat{Z}) = 0$ for all $z \in \widehat{Z}$. \end{itemize} \end{definition}

The following is from Theorem 11 of \cite{LM} and Table 1 of \cite{Bu}.

\begin{proposition}\label{p.Legendre}
A highest weight variety $Z \subset \BP V$ is a homogeneous Legendrian variety if and only if it is isomorphic to one of the following.
The labeling by simple Lie algebras in the  list below is explained in Proposition \ref{p.contactlines} (iii).

\begin{itemize}
    \item[($\fg_2$)] The twisted cubic curve in $\BP^3$.
        \item[($\ff_4$)]  The Lagrangian Grassmannian of 3-dimensional isotropic subspaces in  a symplectic vector space of dimension 6, with the Pl\"ucker embedding.
\item[($\fe_6$)]  The Grassmannian  of 3-dimensional subspaces in a  vector space of dimension 6, with the Pl\"ucker embedding.
    \item[($\fe_7$)] The 15-dimensional highest weight variety  of the $32$-dimensional spin representation of $\fso_{12}.$
\item[($\fe_8$)] The 27-dimensional highest weight variety of the $56$-dimensional irreducible representation of the simple Lie algebra of type  $\textsf{E}_7$.
\item[($\fso_{m}$)] the Segre product $\BP^1 \times \Q^{m-6} \subset \BP^{2m-9}, m \geq 7,$ where $\Q^{m-6} \subset \BP^{m-5}$ is the nonsingular quadric hypersurface of dimension $m-6$.
\end{itemize}
\end{proposition}

We have the following result of Merkulov and Schwachh\"ofer from \cite{MS}.

\begin{theorem}\label{t.MS}
Let $Z \subset \BP V$ be the highest weight variety, which is different from  one of the projective varieties listed in Propositions \ref{p.KN} and  \ref{p.Legendre}.
Let $\sC \subset \BP TM$ be a $Z$-isotrivial cone structure and let $\sP \subset \Fr M$ be the associated G-structure. If  $\sP$ is  torsion-free, then it is locally symmetric.
    \end{theorem}

    As mentioned in the introduction of \cite{HL}, Theorem \ref{t.MS} follows from Theorem 6.12 of \cite{MS}. In fact, the Lie algebra of $\Aut(\widehat{Z})$ is of the form $\C {\rm Id}_V \oplus \fh$ for some semisimple Lie algebra $\fh$. If the structure is not locally symmetric, then $\dim \BK(\fh) >1$ from Lemma 2.3 of \cite{MS}. Thus by Theorem 6.12 of \cite{MS}, the semisimple Lie algebra $\fh$ must be one of those in Table 6 or Table 7 of \cite{MS}. Then one can check that the corresponding highest weight variety is one of those  listed in Propositions \ref{p.KN} and \ref{p.Legendre}.

Combining Theorems \ref{t.Spencer} and \ref{t.MS} with Lemma \ref{l.Bott} (ii),
    we obtain the following.

    \begin{corollary}\label{c.MS}
    Let $Z \subset \BP V$ be a highest weight variety, not one of those  listed in Propositions \ref{p.KN} and \ref{p.Legendre}. Assume that
    \begin{itemize}
    \item[(a)] the homomorphism $\aut(\widehat{Z}) \otimes H^0(Z, \sO(1)) \to H^0(Z, TZ(1))$ is surjective; and
    \item[(b)] $\Xi_Z \subset \partial( \Hom(V, \aut(\widehat{Z})))$ in the notation of Definitions \ref{d.Spencer} and \ref{d.Xi}.
    \end{itemize}
 If a $Z$-isotrivial cone structure $\sC \subset \BP TM$ on a complex manifold $M$ admits  a characteristic conic connection, then the associated G-structure $\sP \subset \Fr M$ is locally symmetric.
    \end{corollary}

To what extent Corollary \ref{c.MS} determines the G-structure?
The local equivalence of a locally symmetric principal connection $\sH$ is determined by the isomorphism type of its curvature tensor $R^{\sH}_x \in \Hom(\wedge^2 T_x M, \sG_x)$ at one point $x \in M$ in the following sense.

\begin{proposition}\label{p.equiv}
Let $\sP \subset \Fr M$ and $\widetilde{\sP} \subset \Fr \widetilde{M}$ be G-structures with the same structure group. Let $\sH$ (resp. $\widetilde{\sH}$) be locally symmetric principal connections on $\sP$ (resp. $\widetilde{\sP})$. Suppose that there exist points $x \in M$ and $\widetilde{x} \in \widetilde{M}$ and a linear isomorphism $\varphi_x: T_x M \to T_{\widetilde{x}} \widetilde{M}$ such that the induced isomorphism $$ \Hom(\wedge^2 T_x M, \End(T_x M)) \stackrel{\varphi_x}{\longrightarrow} \Hom(\wedge^2 T_{\widetilde{x}} \widetilde{M}, \End(T_{\widetilde{x}} \widetilde{M}))$$
sends the curvature $R^{\sH}_x$ at $x$ to the curvature $R^{\widetilde{\sH}}_{\widetilde{x}}$ at $\widetilde{x}$. Then there exists a biholomorphic map $\varphi: O \to \widetilde{O}$ between some neighborhoods $x \in O \subset M$ and $\widetilde{x} \in \widetilde{O} \subset \widetilde{M}$ such that ${\rm d}_x \varphi = \varphi_x$ and the induced map $\varphi_*: \Fr O \to \Fr \widetilde{O}$ sends $\sP|_O$ to $\widetilde{\sP}|_{\widetilde{O}}$. \end{proposition}

Proposition \ref{p.equiv} is an immediate consequence  of two  classical results, Theorem 7.2 in Chapter VI of \cite{KN} saying that the local equivalence type of a locally symmetric affine connection is determined by  the isomorphism type of the curvature tensor at one point and   Example 2.2 in Chapter XI of \cite{KN} saying  that locally symmetric affine connections determine the local structure of symmetric spaces.

By Proposition \ref{p.BK}, the curvature $R^{\sH}_x$ is isomorphic to an element of the space $\BK(\fg)$ of formal curvature maps of the Lie algebra $\fg \subset \End(V)$ of the structure group of the G-structure. Thus the following is a direct corollary of Proposition \ref{p.equiv}.

\begin{corollary}\label{c.equiv}
Let $\fg \subset \End(V)$ be the Lie algebra of a closed subgroup $G \subset {\rm GL}(V)$ such that  $ \dim \BK (\fg) =1$. Let $\sP \subset \Fr M$ and (resp. $\widetilde{\sP} \subset \Fr \widetilde{M}$) be a $G$-structure equipped with a locally symmetric principal connection $\sH$ (resp. $\widetilde{\sH}$) whose curvature $R^{\sH}$ (resp. $R^{\widetilde{\sH}}$) is not zero.
Then for any $x \in M$ and $\widetilde{x} \in \widetilde{M}$,  there exists a biholomorphic map $\varphi: O \to \widetilde{O}$ between some neighborhoods $x \in O \subset M$ and $\widetilde{x} \in \widetilde{O} \subset \widetilde{M}$ such that the induced map $\varphi_*: \Fr O \to \Fr \widetilde{O}$ sends $\sP|_O$ to $\widetilde{\sP}|_{\widetilde{O}}$. \end{corollary}

\begin{definition}\label{d.homogeneous}
An isotrivial cone structure $\sC \subset \BP TM$ on a complex manifold $M$ is
{\em locally flat} (resp. {\em locally symmetric}) if its associated G-structure  is locally flat (resp. locally symmetric). It is
{\em locally homogeneous} if for any two points $x, \widetilde{x} \in M,$ there exist neighborhoods $x \in O, \widetilde{x} \in \widetilde{O}$ and a biholomorphic map $f: O \to \widetilde{O}$ such that its differential ${\rm d} f: \BP TO \to \BP T \widetilde{O}$ sends $\sC|_{O}$ to $\sC|_{\widetilde{O}}$.\end{definition}

We can reformulate Corollary \ref{c.equiv} as follows.

\begin{corollary}\label{c.equivcone}
Let $Z \subset \BP V$ be a submanifold such that $\dim \BK(\aut(\widehat{Z})) =1$. Let $\sC \subset \BP TM$ (resp. $\widetilde{\sC} \subset \BP T\widetilde{M}$) be a locally symmetric $Z$-isotrivial cone structure which is not locally flat. Then for any
$x \in M$ and $ \widetilde{x} \in \widetilde{M},$ there exist neighborhoods $x \in O \subset M, \widetilde{x} \in \widetilde{O} \subset \widetilde{M}$ and a biholomorphic map $f: O \to \widetilde{O}$ such that its differential ${\rm d} f: \BP TO \to \BP T \widetilde{O}$ sends $\sC|_{O}$ to $\widetilde{\sC}|_{\widetilde{O}}.$
\end{corollary}

Another consequence of  Proposition \ref{p.equiv} is the following, which is a special case of Corollary 7.5 in Chapter VI of \cite{KN}.

\begin{corollary}\label{c.homogeneous}
A locally symmetric isotrivial cone structure is locally homogeneous. \end{corollary}

\section{Geometry of adjoint varieties}\label{s.adjoint}

In this section, we use $\fg$ to denote a simple Lie  algebra
whose Dynkin diagram is not of the type $\textsf{A}_1$ or $\textsf{C}_{\ell}$. In other words, we assume that $\fg$ is one of the following types.
$$\textsf{A}_{\ell  \geq 2}, \textsf{B}_{\ell \geq 3}, \textsf{D}_{\ell \geq 4}, \textsf{G}_2, \textsf{F}_4, \textsf{E}_6, \textsf{E}_7, \textsf{E}_8.$$
We use $G$ to denote  the adjoint group of $\fg$, regarded as a subgroup $G \subset {\rm GL}(\fg)$ acting on $\fg$ by the adjoint representation. Thus $\fg$ and $G$ have  different meanings  from the previous sections. This should not cause any confusion.

\begin{definition}\label{d.adjoint}
The {\em adjoint variety} of $\fg$ is a highest weight variety $ Y \subset \BP \fg$ with   $\Aut(\widehat{Y}) = \C^{\times} {\rm Id}_{\fg} \cdot G.$
The underlying complex manifold of an adjoint variety is  called a {\em homogeneous contact manifold} because of the next proposition. \end{definition}

\begin{proposition}\label{p.contact}
For an adjoint variety $Y \subset \BP \fg$,
there is
 a $G$-invariant contact structure on $Y$, namely, a vector subbundle $D \subset TY$  satisfying \begin{equation}\label{e.D} 0 \longrightarrow D \longrightarrow TY \stackrel{\vartheta}{\longrightarrow} \sO(1) \longrightarrow 0, \end{equation} such that
the homomorphism
 $\omega: \wedge^2 D \to \sO(1)$ defined by the Lie brackets of local sections of $D$  gives a nondegenerate anti-symmetric form (i.e. a symplectic form)  $\omega_x : \wedge^2 D_x \to \widehat{x}^*$ for each $x \in Y$. \end{proposition}

\begin{remark}\label{r.typeC}
If $\fg$ has Dynkin diagram of type $\textsf{A}_1$ or $\textsf{C}_{\ell}$, then the highest weight orbit $Y \subset \BP \fg$ is isomorphic to the second Veronese embedding of the projective space $\BP^{\ell}$, the case (III) of Proposition  \ref{p.KN}. So it is not interesting when we try to find applications of Corollary \ref{c.MS}. Moreover,  there are no lines lying on $Y$ and many of the arguments we use in this section cannot work for these cases.    This is why we  exclude them in our discussion. \end{remark}

\begin{example}\label{ex.sl}
When $\fg = \fsl(W)$ for a vector space $W$ of $\dim W \geq 3$,
the adjoint variety has the following description. \begin{itemize}
\item[(i)] The adjoint variety $Y$ of $\fsl(W)$ is biholomorphic to $\BP (T^* \BP W)$, the projectivization of the cotangent bundle of the projective space $\BP W$. \item[(ii)]
We have the natural embedding $$\BP (T^* \BP W) \subset \BP W \times \BP W^*$$
as the correspondence space between  a point in $\BP W$ and a hyperplane in $\BP W^*$, which gives two natural projections $p: Y \to \BP W$ and $q: Y \to \BP W^*$. \item[(iii)]
The direct sum $T^p \oplus T^q \subset TY$ is the contact bundle $D$. Both $T^p_x$ and $T^q_x$ are Lagrangian subspaces of $D_x$ with respect to the symplectic form $\omega_x: \wedge^2 D_x \to \widehat{x}^*$, which induces a perfect pairing $T_x^p \otimes T_x^q \to \widehat{x}^*$.
\item[(iv)] The hyperplane section of the Segre embedding $\BP W \times \BP W^* \subset \BP \End(W)$ by the hyperplane $\fsl(W) \subset \End(W)$ realizes $Y \subset \BP \fsl (W)$ as the adjoint variety of the Lie algebra $\fg = \fsl(W)$. \end{itemize} \end{example}

In the next two propositions, we recall some results on the geometry of lines on adjoint varieties. The following proposition is well-known (e.g. Section 2 of \cite{Hw97}).

\begin{proposition}\label{p.contactlines}
Let $Y$ and $D$ be as in Proposition \ref{p.contact}.
For each $x \in Y$, let $\sC_x \subset \BP T_x Y$ be the closed algebraic subset of tangent vectors to lines through $x$ lying on $Y$.
\begin{itemize}
\item[(i)] All lines on $Y$ are tangent to $D$, namely, the algebraic subset $\sC_x$ is contained in $\BP D_x$ for each $x \in Y$.
    \item[(ii)] The affine cone $\widehat{\sC}_x \subset D_x$ is Lagrangian with respect to the symplectic form $\omega_x$ for each $x \in Y$.
        \item[(iii)] When $\fg$ is $\fsl_m, m \geq 3$, the algebraic subset $\sC_x$ is the disjoint union $\BP T^p_x \cup \BP T^q_x$ described in Example \ref{ex.sl} and for the other simple Lie algebra $\fg$, the algebraic subset $\sC_x \subset \BP D_x$  is the homogenous Legendrian variety  labeled by the Lie algebra $\fg$ in Proposition \ref{p.Legendre}. In particular, the algebraic subset $\sC_x$ is linearly nondegenerate in $\BP D_x,$ namely, the affine cone $\widehat{\sC}_x$ spans $D_x$.
            \item[(iv)] Fix a base point $o \in Y$ and let $P \subset G$ be the isotropy subgroup of $o$. \begin{itemize} \item When $\fg$ is of type $\textsf{A}_{\ell \geq 2}$,   the isotropy representation of $P$ on $\BP D_o$ is
                transitive on $\BP T^p_o$ and $\BP T^q_o$. \item For other $\fg$,   the isotropy representation of $P$ on $D_o$ is irreducible and $\sC_o$ is the highest weight variety. \end{itemize}
  \end{itemize} \end{proposition}

\begin{definition}\label{d.TanLines}
Let $Z \subset \BP V$ be a projective variety. The {\em variety of tangent lines } of $Z$ is the affine cone $\sT_Z \subset \wedge^2 V$ defined by
$$\sT_Z := \{ u \wedge v \mid u \in \check{Z}, v \in T_u \widehat{Z}\}.$$
\end{definition}

The following is Lemma 2.12 of \cite{FL}.

\begin{proposition}\label{p.FuLi}
Let $Y \subset \BP \fg$ be the adjoint variety. Then $\sT_Y$ spans $\wedge^2 \fg$. \end{proposition}

            In the following proposition, (i) and (ii) are from Section 1  of \cite{Hw97} and (iii) is Proposition 4 of \cite{Hw97}. (iv) is a direct consequence of (ii) and (iii).

\begin{proposition}\label{p.TanLines}
Let $Y \subset \BP \fg$ be the adjoint variety.
\begin{itemize}
\item[(i)] The restriction  $TY|_C$ (resp. $D|_C$) to  a line $C \subset Y$ is isomorphic to the direct sum of line bundles $$ \sO(2) \oplus \sO(1)^m \oplus \sO^{m+2} \mbox{ (resp. } \sO(2) \oplus \sO(1)^m \oplus \sO^m \oplus \sO(-1) \mbox{)}$$ on $\BP^1 \cong C$ where $\sO(2) \cong TC$ and $\dim Y = 2m+3$.
    Denote by $D|_C^+ \subset D|_C$ the subbundle corresponding to $\sO(2) \oplus \sO(1)^m$.
    \item[(ii)] For a point $x \in C$ and the corresponding point $\alpha := [T_x C] \in \sC_x$, the affine tangent space $T_{\alpha} \widehat{\sC}_x \subset D_x$ coincides with $(D|_C^+)_x \subset D_x$ in (i).
\item[(iii)] For each $x \in Y$, denote by ${\rm Ker}(\omega_x) \subset \wedge^2 D_x$  the hyperplane annihilated by $\omega_x$.
    Then $\sT_{\sC_x} \subset \wedge^2 D_x$ is contained in ${\rm Ker}(\omega_x)$ and spans it.
    \item[(iv)]
The vector bundle ${\rm Ker}(\omega)$ is spanned by the subset
$$ \bigcup_{C, L}  TC \wedge L \ \subset \ {\rm Ker}(\omega)$$ where we take all lines $C \subset Y$   and all line subbundles $L \subset D|_C^+$ for each line $C \subset Y$. \end{itemize} \end{proposition}

\begin{proposition}\label{p.commute} Let $Y$ and $D$ be as in Proposition \ref{p.contact}. \begin{itemize}
\item[(i)] The composition of natural homomorphisms  $\fg \subset \aut(\widehat{Y}) \to H^0(Y, TY)$ gives an identification $\fg = H^0(Y, TY).$
    \item[(ii)] $H^0(Y, D) = \Hom ( D,  \sO(1)) =0$ and $\Hom(TY, \sO(1)) = \C  \vartheta$.
\item[(iii)] The homomorphism $\vartheta$ in (\ref{e.D}) gives an isomorphism
$H^0(Y, TY) \stackrel{\vartheta}{=} H^0(Y, \sO(1)) = \fg^*.$
\item[(iv)] We have the  following commutative diagram.
$$ \begin{array}{ccc}
H^0(Y, \sO(1)) \otimes H^0(Y, TY) & \longrightarrow & H^0(Y, TY \otimes \sO(1)) \\
\downarrow & & \downarrow \\
H^0(Y, \sO(1)) \otimes H^0(Y, \sO(1)) & \longrightarrow & H^0(Y, \sO(2)),
\end{array} $$
where the horizontal arrows are natural tensor product homomorphisms and the vertical arrows are induced by $\vartheta$.
\item[(v)] We have canonical identifications $$\Hom(D, TY(1)) = \Hom(D, D(1)) = \Hom(D \otimes D, \sO(2)).$$
\end{itemize} \end{proposition}

\begin{proof}
(i) is well-known (e.g. Fact 3.1 of \cite{MS} ). Since $D$ is a contact distribution, the perfect pairing $D \otimes D \to \sO(1)$ induced by $\omega$ gives the isomorphism $$H^0(Y, D) = H^0(Y, D^* \otimes \sO(1)) = \Hom(D, \sO(1))$$ in (ii). If $\phi \in \Hom(D, \sO(1))$,
then the restriction of $\phi: D \to \sO(1)$ to a line $C \subset Y$ induces $\phi|_{TC}: TC = \sO(2)|_C \to \sO(1)|_C$, which must vanish. As tangent vectors to lines span $D$ from Proposition \ref{p.contactlines} (iii), we see $\phi =0$. Thus  $\Hom(D, \sO(1))=0,$ which implies $\Hom(TY, \sO(1)) = \C  \vartheta$.
(iii) follows from (i), (ii), Lemma \ref{l.Bott} and the exact sequence $$0 \longrightarrow H^0(Y, D) \longrightarrow H^0(Y, TY) \longrightarrow H^0(Y, \sO(1)).$$
(iv) is straightforward.
The second identification  in (v) follows from the canonical isomorphism $D^* \otimes  \sO(1) \cong D$ of the perfect pairing induced by $\omega$. To prove the first identification in (v), we need to show that the natural inclusion $\Hom(D, D(1)) \subset \Hom(D, TY(1))$ is an isomorphism, namely,  any $\phi \in \Hom(D, TY(1))$ sends $D$ to into $D(1) \subset TY(1)$.
Choose any line $C \subset Y$. By Proposition \ref{p.TanLines} (i),  $$TY(1)|_C \cong \sO(3) \oplus \sO(2)^m \oplus \sO(1)^{m+2} \mbox{ with } \sO(3) \oplus \sO(2)^m \subset D(1)|_C.$$ Note that $\sO(2)$ in $TY|C$ corresponds to $TC$, which must be sent to $\sO(3) \oplus \sO(2)^m \subset D(1)|_C$ by $\phi$.
As tangent directions to lines span $D$ from Proposition \ref{p.contactlines}, we conclude that $\phi(D) \subset D(1).$
\end{proof}

\begin{proposition}\label{p.wedgeD}
When $\fg$ is not of type $\textsf{A}_{\ell \geq 3}$, the homomorphism $\omega: \wedge^2 D \to \sO(1)$ induces an isomorphism  $$\Hom(\wedge^2 D, \sO(2)) = \Hom(\sO(1), \sO(2)) \cong \fg^*.$$ \end{proposition}

\begin{proof}
Let us first consider the case $\fg$ of type $\textsf{A}_2$. Then $\wedge^2 D$ is a line bundle and $\omega$ is an isomorphism. Thus  $$\Hom(\wedge^2 D, \sO(2)) = \Hom(\sO(1), \sO(2)) \cong \fg^*$$ is immediate.

Assume that $\fg$ is not of type $\textsf{A}_{\ell \geq 1}$.
For any $\phi \in \Hom(\wedge^2 D, \sO(2))$ and  any line $C \subset Y,$ Proposition \ref{p.TanLines} (i) says
$$\phi(TC \wedge D|_C^+) = \phi(\sO(2) \wedge \sO(1)^m) \subset \sO(2),$$ which shows that $\phi (TC \wedge L) =0$ for any line subbundle $L \subset D|_C^+$.  This implies $\phi({\rm Ker}(\omega)) =0$ by Proposition \ref{p.TanLines} (iv). Then $\phi$ factors through an element in $\Hom(\sO(1), \sO(2)) = H^0(Y, \sO(1)) = \fg^*$. It is clear that any element of $\Hom(\sO(1), \sO(2))$ induces an element of $\Hom(\wedge^2 D, \sO(2))$, proving the proposition. \end{proof}

\begin{proposition}\label{p.sym2D}
For each point $x$ of  the adjoint variety $Y \subset \BP \fg$,
let $$Q_x := \{ b \in \Sym^2 D^*_x \mid b ( u, u) = 0 \mbox{ for all } u \in \widehat{\sC}_x \}$$ be the linear subspace of quadrics on $D_x$ vanishing on the algebraic subset $\sC_x \subset \BP D_x.$ Let $Q \subset \wedge^2 D^*$ be the vector subbundle whose fiber at $x$ is $Q_x$. Then
$$\Hom(\Sym^2 D, \sO(2)) = H^0(Y,  \Sym^2 D^* \otimes \sO(2)) = H^0(Y, Q(2)).$$ \end{proposition}

\begin{proof}
We have a natural inclusion $H^0(Y, Q(2)) \subset H^0(Y, \Sym^2 D^* \otimes \sO(2))$. To prove the proposition, it suffices to show that
 any $\phi \in \Hom(\Sym^2 D, \sO(2))$ belongs to $H^0(Y, Q(2))$. But for any line $C \subset Y$,
$$\phi(TC \otimes TC) = \phi(\sO(2) \otimes \sO(2)) \subset \sO(2)$$ shows that $\phi (TC \otimes TC) =0$. It follows that $\phi $ belongs to $H^0(Y, Q(2))$.
\end{proof}

\section{Proof of Theorem \ref{t.adjoint1}}\label{s.Spencer}

Throughout this section, let $\fg$ be as in Section \ref{s.adjoint}, namely, a simple Lie algebra excluding $\textsf{A}_1$ or $\textsf{C}_{\ell}$.
The goal of this section is to prove the following two theorems.

\begin{theorem}\label{t.tensor}
Let $Y \subset \BP \fg$ be the adjoint variety of a simple Lie algebra $\fg$. Then the tensor product homomorphism
$$\Phi: H^0(Y, TY) \otimes H^0(Y, \sO(1)) \to H^0(Y, TY \otimes \sO(1))$$ is surjective. \end{theorem}

\begin{theorem}\label{t.adjointS}
Assume that  $\fg$ is not of type  $\textsf{A}_{\ell \geq 3}.$ For the adjoint variety $Y \subset \BP \fg$, the image of the Spencer homomorphism $\p: \Hom (\fg, \aut(\widehat{Y})) \to \Hom( \wedge^2 \fg, \fg) $ contains $\Xi_Y$. \end{theorem}

Combined with Corollary \ref{c.MS}, they imply the following, proving Theorem \ref{t.adjoint1}.

\begin{corollary}\label{c.adjoint}
Let $Y \subset \BP \fg$ be the adjoint variety of $\fg$, which is not of type $\textsf{A}_{\ell \geq 3}$. Then any $Y$-isotrivial cone structure admitting a characteristic conic connection is locally symmetric. \end{corollary}

It is necessary to exclude $\fg$ of type $\textsf{A}_{\ell \geq 3}$, in Theorem \ref{t.adjointS} and  Corollary \ref{c.adjoint}: counterexamples are discussed in
 Example \ref{ex.A}.

 The proof of Theorem \ref{t.tensor} is easy:

\begin{proof}[Proof of Theorem \ref{t.tensor}]
Let us use the following notation. For $v \in \fg$, let $\vec{v} \in H^0(Y, TY)$ be the corresponding vector field on $Y$ by Proposition \ref{p.commute} (i) and  let $v^* \in H^0(Y, \sO(1)) = \fg^*$ be the corresponding section of $\sO(1)$  by Proposition \ref{p.commute} (iii).

First, consider the case when $\fg$ is of type $\textsf{A}_{\ell \geq 2}$. We have the decomposition $D = T^p \oplus T^q$ from Example \ref{ex.sl}. Fix a base point $o \in Y$.
As $Y$ is homogeneous under $G$, we can choose a general element $u \in \fg$ and  nonzero elements $ v, w \in \fg$ such that \begin{itemize} \item[(a)] $\vec{v}_o \in T^p_o$ and $\vec{w}_o \in T^q_o$;  \item[(b)] $\vec{u}_o \not\in D_o$ and $u^* \in H^0(Y, \sO(1))$ does not vanish at $o \in Y$. \end{itemize}
Set $$v^+ : = \Phi (u^* \otimes \vec{v} - v^* \otimes \vec{u}) \  \in \ H^0(Y, TY (1)),$$ $$w^+ : = \Phi (u^* \otimes \vec{w} - w^* \otimes \vec{u}) \  \in \ H^0(Y, TY (1)).$$ Note that the set of points  $x \in Y$ such that  \begin{itemize} \item the three vectors  $\vec{u}_x$, $\vec{v}_x$ and $\vec{w}_x$ are linearly independent in $T_x Y$ and \item
each of the three sections  $u^*, v^*, w^* \in H^0(Y, \sO(1))$ of the line bundle $\sO(1)$ on $Y$  does not vanish at $x$ \end{itemize} is Zariski-open in $Y$, because $\vec{u}_o, \vec{v}_o$ and $\vec{w}_o$ are linearly independent in $T_o Y$.
It follows that $v^+ \neq 0 \neq w^+ $.
But by Proposition \ref{p.commute} (iii), both $u^* \otimes \vec{v} - v^* \otimes \vec{u}$ and $u^* \otimes \vec{w} - w^* \otimes \vec{u}$ are sent to zero in $H^0(Y, \sO(2))$ under $$0 \longrightarrow H^0(Y, D(1)) \longrightarrow H^0(Y, TY(1)) \longrightarrow H^0(Y, \sO(2)).$$
Thus  $$v^+, w^+ \in H^0(Y, D(1)) = H^0(Y, T^p(1)) \oplus H^0(Y, T^q(1)).$$ Since $u^* \in H^0(Y, \sO(1))$ does not vanish at $o \in Y$, we see that $v^+$ (resp. $w^+$) does not belong to $ H^0(Y, T^q (1))$ (resp. $H^0(Y, T^p(1))$.)  As $$H^0(Y, D(1)) = H^0(Y, T^p(1)) \oplus H^0(Y, T^q (1))$$ is the decomposition into irreducible $G$-modules by Lemma \ref{l.Bott} and  Proposition \ref{p.contactlines} (iv),  the $G$-submodule of $H^0(Y, D(1))$ containing both $v^+$ and  $w^+$ should be the whole $H^0(Y, D(1)).$
Thus  the image of $\Phi$ must contain
$H^0(Y, D(1)).$ Since the image of $\Phi$ is surjective onto $H^0(Y, \sO(2))$ in Proposition \ref{p.commute} (iv), we see that $\Phi$ is surjective.

The proof is simpler when $\fg$ is not of type $\textsf{A}_{\ell}$ because $H^0(Y, D(1))$ is irreducible by Lemma \ref{l.Bott} and Proposition \ref{p.contactlines} (iv). Fix a base point $o\in Y$ and choose $u, v \in \fg$ such that \begin{itemize} \item[(a)] $\vec{v}_o \in D_o$;  \item[(b)] $\vec{u}_o \not\in D_o$ and $u^* \in H^0(Y, \sO(1))$ does not vanish at $o \in Y$. \end{itemize}
Set $v^+ = \Phi(u^* \otimes \vec{v} - v^* \otimes \vec{u})$ as before, then $0 \neq v^+ \in H^0(Y, D(1))$ by the same argument as before. By the irreducibility of $H^0(Y, D(1))$, we see that $H^0(Y, D(1)) \subset {\rm Im}(\Phi)$ and obtain the surjectivity of $\Phi$ from Proposition \ref{p.commute} (iv).
\end{proof}

The proof of Theorem \ref{t.adjointS} consists of several components. The first step is
to define a subspace $S \subset \Hom(\fg, \aut(\widehat{Y}))$ that would give  $\p S = \Xi_Y$.

\begin{definition}\label{d.Sigma}
Let $B$ be the Killing form of $\fg$. \begin{itemize}
\item[(i)] For a bilinear form $b \in \Hom(\fg \otimes \fg, \C) =  \fg^* \otimes \fg^*$ on $\fg$,  define the endomorphism  $b^{\flat} \in \End(\fg)$ by
   $$b (v, w) := B( b^{\flat}(v), w) \mbox{ for all } v, w \in \fg.$$
\item[(ii)] For each $v \in \fg$, let $\ad_v \in \End(\fg)$ be the endomorphism $\ad_v (w) = [v, w]$ for all $w \in \fg$. Denote by $\ad_{\fg} \subset \End(\fg)$ the subspace $\{ \ad_v \mid v \in \fg\}.$ Since $B([v,w], u) + B(w, [v,u]) =0$ for all $u, v, w \in \fg$, we know $\ad_{\fg} \subset (\wedge^2 \fg^*)^{\flat}$.
    \item[(iii)] Let $\Sigma \subset \Sym^2 \fg^*$ be the set of quadrics on $\fg$ vanishing on $Y,$ namely, $$\Sigma := \{ b \in \Sym^2 \fg^* \mid b( v, v) = 0 \mbox{ for all } v \in \widehat{Y}\}.$$  \item[(iv)]
        Define the subspace $S$ of \begin{eqnarray*} \Hom(\fg, \aut(\widehat{Y})) &=& \Hom(\fg, \C {\rm Id}_{\fg} \oplus \fg) \\ &=& \Hom(\fg, \C {\rm Id}_{\fg}) \oplus (\wedge^2 \fg^*)^{\flat} \oplus (\Sym^2 \fg^*)^{\flat} \end{eqnarray*} by $S := \Hom(\fg, \C {\rm Id}_{\fg}) \oplus \ad_{\fg} \oplus \Sigma^{\flat}.$
        \end{itemize} \end{definition}

 \begin{proposition}\label{p.Sigma}
For $S$ as in Definition \ref{d.Sigma}, we have $\p S \subset \Xi_Y$. \end{proposition}

We need the following two lemmata. The first lemma can be seen easily from Section 3.1 and Section 4.2 of \cite{Ya}.

\begin{lemma}\label{l.Killing}  Let $B \in \Sym^2 \fg^*$ be the Killing form of $\fg$. For any $u \in \check{Y} \subset \fg,$
there is a  direct sum decomposition $$\fg = \fg_2 \oplus \fg_1 \oplus \fg_0 \oplus \fg_{-1} \oplus \fg_{-2}$$ as a graded Lie algebra
such that
\begin{itemize}  \item
$\C u = \fg_2$;  \item there exists a nonzero element $E \in \fg_0$ satisfying $$ [\fg_2, \fg_{-2}] = \C E, \ B(E, E) \neq 0 \mbox{ and } [E, v] = i v $$  for  any  $ -2 \leq i \leq 2$ and  $ v \in \fg_i ;$   \item $\fg_0 = \fl \oplus \C E$ for some semisimple Lie algebra $\fl$; \item $B(\fg_i, \fg_j) =0 $ if $i+j \neq 0$ and $B|_{\fg_i \times \fg_{-i}}$ is nondegenerate for each $i = 0, 1, 2.$  \end{itemize} \end{lemma}

\begin{lemma}\label{l.Gu}
         Let $B \in \Sym^2 \fg^*$ be the Killing form of $\fg$. For $u \in \check{Y} \subset \fg$, let $\widehat{D}_u \subset T_u \widehat{Y}$ be the hyperplane corresponding to $D_{[u]} \subset T_{[u]} Y$ and set $$(T_u \widehat{Y})^{\perp} := \{ w \in \fg \mid B(w, T_u \widehat{Y}) =0\}.$$ Then
         \begin{itemize} \item[(i)]  $[\fg, u] =T_u \widehat{Y}$;
         \item[(ii)] $[T_u \widehat{Y}, u] \subset \C u$;
         \item[(iii)]  $[(T_u \widehat{Y})^{\perp}, T_u \widehat{Y}] \subset T_u \widehat{Y}$; and
             \item[(iv)] $[u, \widehat{D}_u] =0.$
         \end{itemize} \end{lemma}

         \begin{proof} (i) is a direct consequence of  the fact  that $Y$ is the  $G$-orbit of  $[u] \in \BP \fg$.
Applying Lemma \ref{l.Killing} to $u$,  we obtain  $$ T_u \widehat{Y} = [\fg, \fg_2] = \fg_2 \oplus \fg_1 \oplus [u, \fg_{-2}]$$ $$(T_u \widehat{Y})^{\perp} = \fg_2 \oplus \fg_1 \oplus (\fg_0 \cap E^{\perp}).$$  Then  $$[T_u \widehat{Y}, u] = [\fg_2 \oplus \fg_1 \oplus [u, \fg_{-2}], \fg_2] \subset \fg_2,$$  proving (ii).
Furthermore, \begin{eqnarray*} [(T_u \widehat{Y})^{\perp}, T_u \widehat{Y}]
&= & [\fg_2 \oplus \fg_1 \oplus (\fg_0 \cap E^{\perp}), \fg_2 \oplus \fg_1 \oplus \C E] \\ & = & \fg_2 \oplus \fg_1 \subset T_u \widehat{Y},\end{eqnarray*} proving (iii).  The contact distribution is defined by the Lie bracket, meaning $$\widehat{D}_u =\{ v \in T_u \widehat{Y} \mid [v, u] =0\}.$$
Thus $\widehat{D}_u = \fg_2 \oplus \fg_1$. Then $$[u, \widehat{D}_u] \in [\fg_2, \fg_2 \oplus \fg_1] =0,$$ proving (iv).
        \end{proof}

         \begin{proof}[Proof of Proposition \ref{p.Sigma}]
         If $h \in \Hom(\fg, \C {\rm Id}_{\fg})$, then $$\p h (u, v) = h(u) \cdot v - h(v) \cdot u \  \in \ \C v + \C u.$$
        Thus $\p h (u, T_u \widehat{Y}) \subset T_u \widehat{Y}$ for any  $u \in \check{Y}$, proving $\p \Hom(\fg, \C {\rm Id}_{\fg}) \subset \Xi_Y.$

        For any $w \in \fg, u \in \check{Y}, v \in T_u \widehat{Y}$,
        $$\p \ad_w (u, v) = [[w, u], v] - [[w,v],u] = [ w, [u, v]] \in [w, \C u] \subset T_u \widehat{Y}$$ by Lemma \ref{l.Gu} (i) and (ii). This shows $\p (\ad_{\fg}) \subset \Xi_Y$.

        For $q \in \Sigma, u \in \check{Y}$ and $ v \in T_u \widehat{Y}$, choose an arc $\{ u + tv + O(t^2) \in \check{Y} \mid t \in \Delta\}$. Then
        $$0 = q(u + t v + \cdots, u + tv + \cdots) = q(u, u) + 2t q(u,v) + O(t^2) $$
        shows $q(u, v) =0$. It follows that $q^{\flat}(u) \in ( T_u \widehat{Y})^{\perp}.$
        Then $$\p q^{\flat}(u, v) = [q^{\flat}(u), v] -[q^{\flat}(v), u] \in [( T_u \widehat{Y})^{\perp}, T_u \widehat{Y}] + [\fg, u] \subset T_u \widehat{Y}$$
        by Lemma \ref{l.Gu} (i) and (iii), proving $\p \Sigma^{\flat} \subset \Xi_Y.$
        \end{proof}

To prove $\p S = \Xi_Y$, we need the following definitions.

\begin{definition}\label{d.zeta}
For a submanifold $Z \subset \BP V$ and a point $z \in \check{Z} \subset V$,
recall that $TZ(-1)_{[z]} = T_z \widehat{Z}/\C z$. Hence $$\Hom(\sO(-1) \otimes TZ(-1), TZ(-1))_{[z]} = \Hom(\C z \otimes T_z \widehat{Z}/ \C z,  T_z \widehat{Z}/ \C z).$$   Define \begin{itemize} \item[(i)]
the homomorphism $\zeta_Z : \Xi_Z  \to \Hom(\sO(-1) \otimes TZ(-1), TZ(-1))$ by sending $\sigma \in \Xi_Z \subset \Hom(\wedge^2 V, V)$ to $$\zeta_Z(\sigma)( z \otimes (v \mod \C z) ) := (\sigma(z, v) \mod \C z) \ \in T_z \widehat{Z}/ \C z$$ for all $z \in \check{Z}$ and $v \in T_z \widehat{Z}$;
\item[(ii)]
the subspace $\Xi'_Z \subset \Xi_Z \subset \Hom(\wedge^2 V, V)$ by $$\Xi'_Z := \{ \sigma \in \Hom(\wedge^2 V, V) \mid \sigma(z, T_z \widehat{Z}) \subset \C z \mbox{ for any } z \in \check{Z}\}$$ such that $\Xi'_Z = {\rm Ker}(\zeta_Z)$; and
\item[(iii)] the homomorphism $\eta_Z: \Xi'_Z \to \Hom(\sO(-1) \otimes TZ(-1), \sO(-1))$ by requiring
$$\eta_Z (\sigma) ( z \otimes (v \mod \C z))  = \sigma(z, v) \in  \C z$$ for all $\sigma \in \Xi'_Z,  z \in \check{Z}$ and $v \in T_z \widehat{Z}$.
 \end{itemize} It is clear that $\eta_Z$ is injective if $\sT_Z$ from Definition \ref{d.TanLines} spans $\wedge^2 V$. \end{definition}

\begin{proposition}\label{p.eta}
Let $Y \subset  \BP \fg$ be the adjoint variety.   Then $\Xi'_Y \subset \Hom(\wedge^2 \fg, \fg)$ is the 1-dimensional vector space generated by the Lie bracket $[,]: \wedge^2 \fg \to \fg$ of the Lie algebra $\fg$ and   $\eta_Y: \Xi'_Y \to \Hom(TY, \sO(1))$ is an isomorphism. \end{proposition}

\begin{proof}
We have $[,] \in \Xi'_Z$ by Lemma \ref{l.Gu} (ii) and $\dim \Hom(TY, \sO(1)) =1$ by Proposition \ref{p.commute} (ii).
Since $\eta_Y$ is injective by Proposition \ref{p.FuLi}, the proposition follows.
\end{proof}

\begin{proposition}\label{p.Psi}
Using Proposition \ref{p.Sigma}, let $\Psi: S \to \Hom(TY, TY(1))$ be the composition
$$\Psi: S \stackrel{\p}{\longrightarrow} \Xi_Y \stackrel{\zeta_Y}{\longrightarrow} \Hom (TY, TY(1)).$$ Consider the following exact sequence arising from (\ref{e.D})
$$0 \longrightarrow \Hom(\sO(1), TY(1)) \longrightarrow \Hom(TY, TY(1)) \stackrel{r_D}{\longrightarrow} \Hom(D, TY(1)).$$
Then \begin{itemize}
\item[(i)] ${\rm Ker}(\Psi) = \C B^{\flat} \subset \Sigma^{\flat}$ and $\p (\C B^{\flat}) = \Xi'_Y = {\rm Ker}(\zeta_Y).$
    \item[(ii)] $\Psi$ sends $\ad_{\fg}$ isomorphically onto $\Hom(\sO(1), TY(1))$.
  \end{itemize}
\end{proposition}

\begin{proof}
If $B \not\in \Sigma$, the quadric hypersurface defined by $B$ in $\BP \fg$ cuts $Y \subset \BP \fg$ along a hypersurface in $Y$ preserved by the $G$-action, a contradiction. Thus $B \in \Sigma$.
Note that $B^{\flat} \in \Hom(\fg, \fg)$ is just ${\rm Id}_{\fg}$. Thus for any $u \in \check{Y}$ and $v \in T_u \widehat{Y}$,
$$\p B^{\flat} (u, v) = [u, v] - [v, u] = 2 \, [u, v] \ \in \C u$$
by Lemma \ref{l.Gu} (ii), verifying $B^{\flat} \in \Xi'_Y$. Since $\p$ is injective from Proposition \ref{p.KN}, the rest of (i) follows from Proposition \ref{p.eta}.

To prove (ii), we claim that $r_D \circ \Psi(\ad_{\fg}) =0.$ In fact,
for any $w \in \fg$ and $u \in \check{Y}$, the homomorphism $r_D(\Psi(\ad_w)) \in \Hom(D, TY(1))$ sends $\widehat{D}_u$ to $[w, [u, \widehat{D}_u]]$, which is zero by Lemma \ref{l.Gu} (iv).
From the claim, we have $\Psi(\ad_{\fg}) \subset \Hom(\sO(1), TY(1))$.
But ${\rm Ker}(\Psi) \cap \ad_{\fg} =0$ by (i) and $\Hom(\sO(1), TY(1)) = H^0(Y, TY) = \fg$. This proves (ii).
 \end{proof}

The next proposition completes the proof of Theorem \ref{t.adjointS}.

\begin{proposition}\label{p.sym}
The homomorphism in  Proposition \ref{p.Psi} induces a homomorphism
$$\Psi':  \Hom(\fg, \C \, {\rm Id}_{\fg}) \oplus (\Sigma^{\flat}/\C B^{\flat}) \longrightarrow
\Hom(D, TY(1)).$$  Assume that $\fg$ is not of type $\textsf{A}_{\ell \geq 3}$.  Then $\Psi'$ is an isomorphism, which, combined with Propositions \ref{p.Sigma} and \ref{p.Psi}, shows that  the image  $\p S$ is equal to $\Xi_Y$
 \end{proposition}

\begin{proof}
By Proposition \ref{p.Psi}, the homomorphism $\Psi'$ is injective. Recall from Proposition \ref{p.commute} (v),
\begin{eqnarray*} \Hom(D, TY(1)) &=& \Hom(D\otimes D, \sO(2)) \\ &=& \Hom(\wedge^2 D, \sO(2)) \oplus \Hom(\Sym^2 D, \sO(2)). \end{eqnarray*} From the assumption that $\fg$ is not of type $\textsf{A}_{\ell}$, we have $$ \Hom(\fg, \C \, {\rm Id}_{\fg}) \cong \fg^* \cong \Hom(\wedge^2 D, \sO(2))$$ by Proposition \ref{p.wedgeD}. Thus to prove that $\Psi'$ is an isomorphism, it suffices to show that $(\Sigma^{\flat}/\C B^{\flat})$ and $\Hom(\Sym^2 D, \sO(2))$ have the same number of irreducible components when they are decomposed into irreducible $\fg$-modules. The number $s$ of $\fg$-irreducible components in  $(\Sigma^{\flat}/\C  B^{\flat})$ is given in Lemma \ref{l.s} below.
Recall $\Hom(\Sym^2 D, \sO(2)) = H^0(Y, Q(2))$ from Proposition \ref{p.sym2D}. Pick a point $x \in Y$ and let $P \subset G$ be the isotropy subgroup at $x$. It is easy to see that the unipotent radical of $P$ acts trivially on $D_x$, thus the $P$-representation on $Q_x$ is completely reducible.  By Lemma \ref{l.Bott}, the number of $\fg$-irreducible components in $H^0(Y, Q(2))$ is equal to the number $s'$ of $P$-irreducible components in $Q_x$, which is given in  Lemma \ref{l.s'} below. Since $s=s'$ from the two lemmata, we conclude that $\Psi'$ is an isomorphism.
\end{proof}

\begin{lemma}\label{l.s}
Let $s$ be the number of irreducible components in the decomposition of $\Sigma/\C B$ as a $\fg$-module. Then $s$ is as follows depending on the type of $\fg$.
\begin{itemize}\item[(i)] $s=2$ for $\textsf{B}_{\ell \geq 3}$, $\textsf{D}_{\ell \geq 5}.$
\item[(ii)] $s=3$ for $\textsf{D}_4$.
\item[(iii)] $s=1$ for $\textsf{A}_2, \textsf{G}_2, \textsf{F}_4, \textsf{E}_6, \textsf{E}_7$ or $\textsf{E}_8$.
   \end{itemize} \end{lemma}

   \begin{proof}
Explicit expressions of the decomposition of $\Sym^2 \fg^*$ into irreducible $\fg$-modules can be found in Table 5 in pages 300 -- 305 of \cite{OV}. We list them below using the notation of \cite{OV} ($\pi_i$ is a fundamental weight,  $1$ denotes the 1-dimensional representation, which corresponds to $\C B$ and  ${\rm Ad}$ denotes the adjoint representation on $\fg$).
\begin{itemize}
\item $\textsf{A}_2:  R(2 \pi_1 + 2 \pi_2) + {\rm Ad} + 1$
      \item $\textsf{B}_3: R(2 \pi_3) + R( 2 \pi_1) + R(2 \pi_2) +1$
    \item $\textsf{B}_4: R(2 \pi_4) + R(2 \pi_1) + R( 2 \pi_2) + 1$
    \item $\textsf{B}_{\ell \geq 5}: R(\pi_4) + R(2 \pi_1) + R(2 \pi_2) + 1$
          \item $\textsf{D}_4:  R(2 \pi_1) + R( 2 \pi_2) + R( 2 \pi_3) + R(2 \pi_4) + 1$
        \item $\textsf{D}_{\ell \geq 5}: R(\pi_4) + R(2 \pi_1) + R(2 \pi_2) + 1$
                          \item $\textsf{G}_2: R(2 \pi_2) + R(2 \pi_1) +1$
                \item $\textsf{F}_4: R(2 \pi_4) + R(2 \pi_1) + 1$
                \item $\textsf{E}_6: R( 2 \pi_6) + R(\pi_1 + \pi_5) +1$
                \item $\textsf{E}_7: T( 2 \pi_6) + R(\pi_2) +1$
                \item $ \textsf{E}_8: R( 2 \pi_1) + R(\pi_7) + 1$
                \end{itemize}
Consider the restriction homomorphism $\Sym^2 \fg^* \to H^0(Z, \sO(2))$ the kernel of which must be $\Sigma$. Since $H^0(Z, \sO(2))$ is an irreducible $\fg$-module by Lemma \ref{l.Bott}, we see that $\Sigma \subset \fg^*$ must contain all irreducible components of the above decomposition, excepting exactly one irreducible component. Thus the value of $s$ must be as stated. \end{proof}

\begin{lemma}\label{l.s'}
Let $x \in Y$ be a point and let $Q_x \subset \Sym^2 D^*_x$ be as in Proposition \ref{p.sym2D}.  Let $s'$ be the be the number of irreducible components in the decomposition of $Q_x$ as a $P$-module. Then $s'$ is as follows depending on the type of $\fg$.
\begin{itemize}
\item[(i)] $s'=2$ for   $\textsf{B}_{\ell \geq 3}$ or $\textsf{D}_{\ell \geq 5}.$
\item[(ii)] $s'=3$ for $\textsf{D}_4$.
\item[(iii)] $s'=1$ for $\textsf{A}_2, \textsf{G}_2, \textsf{F}_4, \textsf{E}_6, \textsf{E}_7$ or $\textsf{E}_8$.
   \end{itemize} \end{lemma}

\begin{proof}
Consider first the case when $\fg$ is of type $\textsf{A}_2$. Then $\sC_x$ is two distinct points on $\BP D_x \cong \BP^1$. Thus $\dim Q_x =1$ and $s'=1$.

When $\fg$ is not of type $\textsf{A}_{\ell}$,  the algebraic subset $\sC_x$ is irreducible and  the  affine cone $\widehat{\sC}_x \subset D_x$ is Lagrangian with respect to the symplectic form $\omega_x$. From Lemma 5.6 of  \cite{Bu}, this implies that   the vector space $Q_x$ is naturally isomorphic to the Lie algebra $\aut(\widehat{\sC}_x)/\C \, {\rm Id}_{D_x}$.
From the list of $\sC_x$ in  Table 1 of  \cite{Bu} or Section 1.4.6 of \cite{Hw01}, we have the following list of the pair  $(\fg : \aut(\widehat{\sC}_o)/\C {\rm Id}_{D_o}).$
\begin{itemize}
\item $\textsf{B}_{\ell}: \textsf{A}_1 \times \textsf{B}_{\ell-2}$ for $\ell \geq 3$
\item $\textsf{D}_4: \textsf{A}_1 \times \textsf{A}_1 \times \textsf{A}_1$
\item $\textsf{D}_{\ell}: \textsf{A}_1 \times \textsf{D}_{\ell -2}$ for $\ell \geq 5$
\item $\textsf{G}_2: \textsf{A}_1$
\item $\textsf{F}_4: \textsf{C}_3$
\item $\textsf{E}_6: \textsf{A}_5$
\item $\textsf{E}_7: \textsf{D}_6$
\item $\textsf{E}_8: \textsf{E}_7$ \end{itemize}
From this list, we obtain the value $s'$ as stated in the lemma.
\end{proof}

\section{Space of formal curvature maps of the adjoint variety}\label{s.Bianchi}

In this section we give a proof of Theorem \ref{t.Bianchi}. The arguments are rather different depending on whether the simple Lie algebra $\fg$ is of type $\textsf{A}_2$ or not. We treat these two cases separately in the two subsections below.

\subsection{Theorem \ref{t.Bianchi} for $\fg$ different from $\textsf{A}_{\ell}$ or $\textsf{C}_{\ell}$}\label{ss.Bianchi}

The goal of this subsection is to prove the following which is exactly Theorem \ref{t.Bianchi} for a simple Lie algebra $\fg$  of type different from $\textsf{A}_{\ell}$ or $\textsf{C}_{\ell}.$

\begin{theorem}\label{t.Bianchi1}
Let $\fg$ be a simple Lie algebra of type different from $\textsf{A}_{\ell \geq 1}$ or $\textsf{C}_{\ell}$. Let $\ad_{\fg} \subset \End(\fg)$ be the image of the adjoint representation and denote by $\widehat{\fg} \subset \End(\fg)$  the direct product Lie algebra   $$ \widehat{\fg} :=  \ad_{\fg} \oplus \C\, {\rm Id}_{\fg}  \ \subset \ \End(\fg).$$ Then   $ \BK(\widehat{\fg})$ is 1-dimensional and  generated by the Lie bracket $$[\,,\,] \in \Hom(\wedge^2 \fg, \ad_{\fg}) \subset \Hom(\wedge^2 \fg, \widehat{\fg}).$$
In other words, if $h \in \Hom(\wedge^2 \fg, \widehat{\fg})$  satisfies $$h^{\sharp}(x,y,z) := h(x,y) \cdot z + h(y,z) \cdot x + h(z,x)\cdot y =0 $$  for all $x,y,z \in \fg$, then $h = s \, [\, ,\,] \in \Hom (\wedge^2 \fg, \ad_{\fg})$ for some $s \in \C.$ \end{theorem}

The proof of Theorem \ref{t.Bianchi1} uses the following fact  on $\End(\fg)$.

\begin{lemma}\label{l.decompo}
There are  irreducible $\fg$-submodules $R, \widetilde{\fg} \subset \wedge^2 \fg$ such that
  \begin{itemize} \item[(i)] $\wedge^2 \fg = R \oplus \widetilde{\fg}$; \item[(ii)] the Lie bracket $[\, ,\,]: \wedge^2 \fg \to \fg$ annihilates $R$ and sends   $\widetilde{\fg}$  isomorphically to $\fg$; and \item[(iii)]  $\widetilde{\fg}  \subset \wedge^2 \fg \subset \End(\fg)$ is the only $\fg$-submodule of $\End(\fg)$ isomorphic to $\fg$. \end{itemize}
\end{lemma}

\begin{proof}
We list the decomposition of $\wedge^2 \fg$  into irreducible $\fg$-modules below, which can be found, for example, Table 5 in pages 300 -- 305 of \cite{OV}. Here we use the notation of \cite{OV},  where $\pi_i$ denotes a fundamental weight, $R(\lambda)$ denotes the irreducible representation with the highest weight $\lambda$ and  ${\rm Ad}$ denotes the adjoint representation on $\fg$.
\begin{itemize}
%\item $\textsf{A}_2: R(3 \pi_1) + R(3 \pi_2) + {\rm Ad}$
 \item $\textsf{B}_3: R(\pi_1 + 2 \pi_3) + {\rm Ad}$
       \item $\textsf{B}_{\ell \geq 4}: R(\pi_1+ \pi_3) + {\rm Ad}$
          \item $\textsf{D}_4:  R(\pi_1 + \pi_3 + \pi_4) + {\rm Ad}$
        \item $\textsf{D}_{\ell \geq 5}:  R(\pi_1+ \pi_3) + {\rm Ad}$
                          \item $\textsf{G}_2: R(3\pi_1) + {\rm Ad}$
                \item $\textsf{F}_4: R(\pi_3) + {\rm Ad}$
                \item $\textsf{E}_6: R( \pi_3) +{\rm Ad}$
                \item $\textsf{E}_7: T( \pi_5) + {\rm Ad}$
                \item $ \textsf{E}_8: R(  \pi_2) + {\rm Ad}$
                \end{itemize}
(i) and (ii) follow from the above description. (iii) follows from the above list combined with the decomposition of $\Sym^2 \fg$ into irreducible $\fg$-submodules recalled in the proof of Lemma \ref{l.s}.
\end{proof}

The crucial step in the proof of Theorem \ref{t.Bianchi1}  is the following proposition.

\begin{proposition}\label{p.hR}
Fix a Borel subalgebra $\fb \subset \fg$.
Let $S \subset \BK(\widehat{\fg})$ be a $\fg$-irreducible submodule and let $0 \neq h \in S$ be a highest weight element with respect to the choice of $\fb$.  Then $h(R) =0$ for $R$ as in Lemma \ref{l.decompo}.  \end{proposition}

In fact,  Theorem \ref{t.Bianchi1} can be deduced from Proposition \ref{p.hR} as follows.

\begin{proof}[Proof of Theorem \ref{t.Bianchi1} assuming Proposition \ref{p.hR}]
For any choice of $S$ as in Proposition \ref{p.hR} and any element $h \in S$,  we claim that $h = s \, [\, , \,]$ for some $s \in \C$. As this holds for any choice of $S$, we obtain $S = \BK (\widehat{\fg}) = \C \,  [\, , \,]$, proving the theorem.

It remains to prove the claim. By Proposition \ref{p.hR}, for any highest weight element $h$ of $S,$ we have $h(R) =0$. But the highest weight elements under different choices of Borel subalgebras span $S$. Thus we can assume that $h(R) =0$ for any $h \in S$. This means that for any $h \in S$,    there exist $\phi\in \Hom(\fg, \ad_{\fg})= \End(\fg)$ and $\eta\in\fg^*$ such that $h =\phi \circ [\, , \, ]+\eta([ \, , \,]) {\rm Id}_{\fg}$. Hence for all $x, y, z\in\fg$, we have
\begin{eqnarray}\label{e.proofBianchi}
0 & = & [\phi([x, y]), z]+[\phi([y, z]), x]+[\phi([z, x]), y] \\
&& +\eta([x, y]) z+\eta([y, z]) x+\eta([z, x]) y. \nonumber
\end{eqnarray}

Choose any $x \in \check{Y}$ in the highest weight variety of $Y \subset \BP \fg$. By Lemma \ref{l.Killing}, we have a gradation
$$\fg=\fg_2\oplus\fg_1\oplus \fg_0 \oplus\fg_{-1}\oplus\fg_{-2} \mbox{ with } \fg_0 = \fl \oplus \C E$$ such that $x \in \fg_2$ and  $[\fl, \fg_2] =0.$    Choose any  $y, z\in\ker\ad_x=\fg_2\oplus\fg_1\oplus\fl$. Then $[\phi([y, z]), x]=-\eta([y, z])x\in\C x$  from (\ref{e.proofBianchi}). It follows that
$$
\phi([y, z])\in\fg_2\oplus\fg_1\oplus\fl\oplus\C E \mbox{ for any } y, z\in\ker\ad_{x}=\fg_2\oplus\fg_1\oplus\fl. $$
Choosing $y, z\in\fg_1$ with $[y, z]=x$,  this gives \begin{equation}\label{e.phiyz} \phi(x)\in \fg_2\oplus\fg_1\oplus\fl\oplus\C E.\end{equation}
Taking $x \in \fg_2, y\in\fl$ and $z=E$ in (\ref{e.proofBianchi}),  we have $0= [\phi(2x), y]+\eta(2x)y$. Hence $[\phi(x), \fl]\subset\fl$, which implies that
\begin{equation}\label{e.phix} \phi(x)\in\fg_2\oplus\C E\oplus\fg_{-2}. \end{equation}
Combining (\ref{e.phiyz}) and (\ref{e.phix}), we obtain $\phi(x) \in \C x \oplus \C E.$
Since $[\fg_2 + \fg_1, \fg_2] =0$ and $[\fg_2 + \fg_1, E] = \fg_2 + \fg_1$, the action of the subgroup $\exp(\fg_2 + \fg_1) \subset {\rm GL}(\fg)$ on $\BP \fg$ fixes the point $[x] \in Y$ and moves $[E] \in \BP \fg$ freely. It follows that we can choose a new gradation
$$\fg=\fg'_2\oplus\fg'_1\oplus (\fl'\oplus\C E')\oplus\fg'_{-1}\oplus\fg'_{-2}$$
with $\fg_2 = \fg'_2$ and $E' \not\in \C E \oplus \fg_2$. Applying the same argument to this new gradation, we obtain  $\phi(x) \in \C x \oplus \C E'.$ It follows that $\phi(x) \in \C x$. Since this holds for all $x$ in the highest weight variety $Y$, we see that $\phi = s \, {\rm Id}_{\fg}$ for some $s \in \C$. Then the first line of (\ref{e.proofBianchi}) is zero by Jacobi identity for the Lie bracket. Thus $$\eta([x, y]) z+\eta([y, z]) x+\eta([z, x]) y =0$$ for all $x,y,z \in \fg$. Putting in nonzero elements $x, y \in \fg_1$ and $z = [x,y] \in \fg_2,$ we have $\eta(z) =0$ for any $z \in \check{Y}$. This implies $\eta =0$, proving the claim.
\end{proof}

The rest of this subsection is devoted to the proof of Proposition \ref{p.hR}.
We need to recall the following from  Sections 3.1 and  4.2 of \cite{Ya}. This is a refined version of Lemma \ref{l.Killing}.

\begin{lemma}\label{l.alphak}
Fix a Cartan subalgebra and a Borel subalgebra $\ft \subset \fb \subset \fg$ and
let $\Phi \subset \ft^*$ be the corresponding root system of $\fg$ with positive roots $\Phi^+$ and simple roots $\alpha_1, \ldots, \alpha_{\ell} \in \Phi^+$. Let $\theta \in \Phi$ be the highest root and let $\fg_{\alpha} \subset \fg$ be the root subspace of a root $\alpha \in \Phi$. For two roots $\alpha, \beta$, we have the integer $\langle \alpha, \beta \rangle := 2 B(\alpha, \beta)/B(\beta, \beta)$ in terms of the Killing form $B \in \Sym^2 \fg^*$. \begin{itemize}
\item[(i)] The highest root $\theta$ is a fundamental weight, more precisely, there exists a long simple root $\alpha_k$ such that $\langle \alpha_j, \theta\rangle = \delta_{jk}$ for all $1 \leq j \leq \ell$.
    \item[(ii)] Set $\Phi_i := \{ \alpha \in \Phi \mid \langle \alpha, \theta \rangle =i\}.$
Then $\Phi$ is the disjoint union of $\Phi_2, \Phi_1, \Phi_0, \Phi_{-1}, \Phi_{-2}$ with $\Phi_2 = \{ \theta \}$ and $\Phi_{-2} =\{ -\theta \}.$
\item[(iii)] Define $\fg_i := \oplus_{\alpha \in \Phi_i} \fg_{\alpha}$ for $i \neq 0$ and $\fg_0 := \ft \oplus (\oplus_{\alpha \in \Phi_0} \fg_{\alpha}).$ Then $$
\fg = \fg_2 \oplus \fg_{-1} \oplus \fg_0 \oplus \fg_{-1} \oplus \fg_{-2}$$ is a graded Lie algebra such that the Lie bracket $[\, , \, ]: \fg_1 \times \fg_1 \to \fg_2$ is a nondegenerate pairing.
\item[(iv)] Write a root $\alpha \in \Phi$ as $\alpha = \sum_{j=1}^{\ell} c_j(\alpha) \, \alpha_j$ in terms of simple roots. Then $c_k(\alpha) = \langle \alpha, \theta \rangle$ and $\Phi_i := \{ \alpha \in \Phi \mid c_k(\alpha) = i\}. $
    \item[(v)] The homomorphisms $\fg_0 \to \End(\fg_1)$ and $\fg_2 \to \Hom(\fg_{-1}, \fg_1)$ induced by the Lie bracket  are injective.  Thus $\theta - \alpha \in \Phi_1$ for any $\alpha \in \Phi_{1}$.
    \end{itemize} \end{lemma}

\begin{remark}\label{r.alphak} The root $\alpha_k$  is given by the node connected to $-\theta$ in the extended Dynkin diagrams in p. 454 of \cite{Ya}. \end{remark}

To proceed,  we need to introduce some notation.

\begin{notation}
Using Lemma \ref{l.alphak},
we choose some root vectors of $\fg$ in the following way.
\begin{itemize}
\item Choose  nonzero vectors $v_{\theta} \in \fg_2$ and
$v_{\alpha_k} \in \fg_{\alpha_k}$.  \item Choose
$v_{-\theta} \in \fg_{-\theta}$ such that $t_{\theta} := [v_{\theta}, v_{-\theta}]$ satisfies $[t_{\theta}, v_{\theta}] = 2 \, v_{\theta}$.
\item Choose $v_{-\alpha_k} \in \fg_{- \alpha_k}$  such that $t_{\alpha_k} := [v_{\alpha_k}, v_{-\alpha_k}]$ satisfies $[t_{\alpha_k}, v_{\alpha_k}] = 2 \, v_{\alpha_k}.$ \end{itemize}
    Note that $\theta- \alpha_k$ and $-\theta + \alpha_k$ are  roots by Lemma \ref{l.alphak} (v). \begin{itemize}
    \item Choose $v_{\theta- \alpha_k} \in \fg_{\theta - \alpha_k}$ such that $[v_{\alpha_k}, v_{\theta-\alpha_k}] = v_{\theta}$.
        \item Choose $v_{-\theta + \alpha_k} \in \fg_{-\theta + \alpha_k}$ such that $[v_{-\alpha_k}, v_{-\theta + \alpha_k}] = v_{-\theta}.$
            \end{itemize}
\end{notation}

For the proof of Proposition \ref{p.hR}, we  need several lemmata.

\begin{lemma}\label{l.wedgefg}
The irreducible $\fg$-module $R \subset \wedge^2 \fg$ in Lemma \ref{l.decompo} has a highest weight vector $v_{\theta} \wedge v_{\theta - \alpha_k}$ and a lowest weight vector $v_{-\theta} \wedge v_{- \theta + \alpha_k}.$ \end{lemma}

\begin{proof}
 From Lemma \ref{l.alphak}, the vector $v_{\theta} \wedge v_{\theta - \alpha_k}$ belongs to $R= {\rm Ker}([\, ,\,])$ and is annihilated by the action of any positive root vector. Thus it is a highest weight vector and $v_{-\theta} \wedge v_{- \theta + \alpha_k}$ is a lowest weight vector.  \end{proof}

\begin{lemma}\label{l.gamma}
 Let $\beta \in \Phi$ be a negative root such that $c_k(\beta) =0$ and $- \beta$ is not a simple root. Then there exists a positive root $\gamma$   such that \begin{itemize} \item[(i)] $c_k(\gamma) =0$; \item[(ii)]   $\beta + \gamma$ is a negative root; and \item[(iii)] $ \langle \beta + \gamma, \alpha_k \rangle \leq 0.$ \end{itemize}
\end{lemma}

\begin{proof}
Since $-\beta$ is not a simple root, there is some simple root $\alpha_j$ such that $\beta':= \beta + \alpha_j$ is a negative root. From $c_k(\beta) =0$, we have $\alpha_j \neq \alpha_k$ and $\langle \alpha_k, \alpha_j \rangle \leq 0$.    We consider the following  two cases separately.

If $\langle \alpha_k, \alpha_j\rangle =0$, set $\gamma := - \beta' = -\beta - \alpha_j$ with $c_k(\gamma) = -c_k(\beta) =0$. Then $\beta + \gamma  =- \alpha_j $ is a negative root and $$\langle \beta + \gamma, \alpha_k\rangle = - \langle \alpha_j, \alpha_k\rangle =0.$$

If $\langle \alpha_k, \alpha_j\rangle \leq -1$, set $\gamma := \alpha_j$ with $c_k(\gamma) =0$. Then $\beta+ \gamma = \beta'$ is a negative root.  Note that $\langle \beta, \alpha_k \rangle \leq 1$ because $\alpha_k$ is a long root.  Thus
 $$\langle \beta + \gamma, \alpha_k\rangle = \langle \beta, \alpha_k\rangle + \langle \alpha_j, \alpha_k \rangle \leq 0.$$
\end{proof}

\begin{lemma}\label{l.lambda}
Let $\lambda$ be the weight of $h \in S$ in Proposition \ref{p.hR}.
\begin{itemize}
\item[(i)] Write $h = f + \sigma \, {\rm Id}_{\fg}$ with $f \in \Hom (\wedge^2 \fg, \fg)= \Hom(\wedge^2 \fg, \ad_{\fg})$ and $\sigma \in \wedge^2 \fg^*$. If $f$ (resp. $\sigma$) is nonzero, then its weight is $\lambda$.
    \item[(ii)] Let $G \subset {\rm GL}(\fg)$ be the adjoint group with Lie algebra $\ad_{\fg}$. Let $B \subset G$ be the Borel subgroup corresponding to $\fb \subset \fg$ and ${\rm Unip}(B) \subset B$ be its unipotent radical with Lie algebra $\fu= \oplus_{\alpha \in \Phi^+} \fg_{\alpha}.$ Then for any $a \in {\rm Unip}(B)$ and $x, y \in \fg$,
        $$a \cdot (h(x,y)) = h(a \cdot x, a \cdot y),  \ a\cdot (f(x,y)) = f(a \cdot  x, a \cdot y), \ \sigma(x, y) =  \sigma(a \cdot x, a \cdot y).$$
        \item[(iii)]  Write $$w:=h(v_{-\theta}\wedge v_{-\theta+\alpha_k}), \ u:=f(v_{-\theta}\wedge v_{-\theta+\alpha_k}) \ \in \fg$$ and $c :=\sigma(v_{-\theta}\wedge v_{-\theta+\alpha_k}) \in \C.$ Then  $w=0$ (resp. $u=0$, resp. $c =0$) if and only if $h(R)=0$ (resp. $f(R)=0$, resp. $\sigma(R)=0$).
\item[(iv)] If $w\neq 0$, it is a $\ft$-eigenvector of weight $\wt(w) =\lambda-2\theta+\alpha_k\in\Phi \cup\{0\}.$
    \item[(v)] If $c \neq 0$, then $\lambda - 2 \theta + \alpha_k =0$.
        \end{itemize} \end{lemma}

\begin{proof}
(i) is obvious from $\widehat{\fg} =\ad_{\fg} \oplus\C{\rm Id}_{\fg}.$

Since $h$ is of highest weight, we have $b \cdot h \in \C h$ for any $b \in B$ and $a \cdot h= h$ for any $a \in {\rm Unip}(B)$. Thus for all $x, y \in \fg$ and $a \in {\rm Unip}(B)$,  $$h(x,y) = (a \cdot h) (x,y) = a \cdot (h(a^{-1} \cdot x, a^{-1} \cdot y))$$  $$f(x,y) = (a \cdot f) (x,y) = a \cdot (f(a^{-1} \cdot x, a^{-1} \cdot y))$$ $$\sigma(x,y) \, {\rm Id}_{\fg} = (a \cdot \sigma) (x, y) \, {\rm Id}_{\fg}  = \sigma(a^{-1} \cdot x, a^{-1} \cdot y) \, {\rm Id}_{\fg},$$ from which (ii) follows.

Since $r:= v_{-\theta}\wedge v_{-\theta+\alpha_k}$ is a  vector in $R$ with the lowest weight,  the orbit  $G\cdot [r] \in \BP R$ is the highest weight variety  spanning $\BP R$. For the Lie algebra $\fu$ of ${\rm Unip}(B)$, we have $$[\fu, r] + \C r = [\fb, r] = [\fg, r]$$ because the subgroup $U^- \subset G$ with Lie algebra $\fu^- := \oplus_{\alpha \in \Phi^-} \fg_{\alpha}$ fixes $r$. Thus the orbit ${\rm Unip}(B) \cdot [r]$ is open in $G \cdot [r]$.  If $w=h(r) =0,$ then $$0 = a \cdot (h( r)) = h ( a \cdot r) \mbox{ for all } a \in {\rm Unip}(B)$$ by (ii). Since ${\rm Unip}(B) \cdot [r]$ spans $ \BP R$, we see that $h(R) =0$.
The conclusions (iii) for $f$ and $\sigma$ follow from the same argument.

(iv) is straightforward. In (v), the weight of $c \neq 0$ must be $\lambda - 2 \theta + \alpha_k$ as in (iv). Since the weight of the constant $c$ must be 0, we have $\lambda - 2 \theta + \alpha_k =0$.
\end{proof}

\begin{lemma}\label{l.exclude}
If $h(R)\neq 0$ in the setting of Lemma \ref{l.lambda}, then $\lambda$ is distinct from $\theta$, $\theta-\alpha_k$ and $2\theta-\alpha_k-\alpha_j$ for any $j=1,\ldots,\ell$.
\end{lemma}

\begin{proof}
From the assumption $h(R) \neq 0$, we have
 $\beta: = \lambda-2\theta+\alpha_k \in \Phi \cup \{ 0 \}$ by Lemma \ref{l.lambda} (iv).
If $\lambda$ is equal to $\theta$ (resp. $\theta-\alpha_k$, resp. $2\theta-\alpha_k-\alpha_j$ for some $j$),  then $\beta$ is equal to $-\theta+\alpha_k$ (resp. $-\theta$, resp. $-\alpha_j$ for some $j$), which implies that  $\beta$ is a negative root. We write  $\delta=-\beta$ and derive a contradiction for each of the following cases for the positive root $\delta$.

\medskip
{\bf Case when $\delta = \theta$}.
Since the weight $\wt(h)=2\theta-\alpha_k-\delta$ is dominant,
$$0\leq\langle\wt(h), \alpha_k\rangle=-\langle\delta,\alpha_k\rangle = - \langle \theta, \alpha_k \rangle = -1,$$
a contradiction.

\medskip
{\bf Case when  $\delta=\alpha_k$}. Since the weight $\wt(h)=2\theta-\alpha_k-\delta$ is dominant, $$0 \leq \langle\wt(h), \alpha_k\rangle= \langle 2 \theta - 2 \alpha_k, \alpha_k \rangle = -2,$$ a contradiction.

\medskip
{\bf Case when $\delta = \alpha_j $ for some  $j\neq k$.} We have
$$0\leq \langle\wt(h), \alpha_j\rangle=\langle 2\theta, \alpha_j\rangle-\langle\alpha_k, \alpha_j\rangle-\langle \alpha_j, \alpha_j\rangle,$$
implying $\langle\alpha_j,\alpha_j\rangle\leq-\langle\alpha_k, \alpha_j\rangle$. Thus $\alpha_j$ is a strictly short root adjacent to $\alpha_k$ in the Dynkin diagram of $\fg$. From  the location on $\alpha_k$ in  the extended Dynkin diagram in page 454 of \cite{Ya} (see Remark \ref{r.alphak}), this is possible  only if $\fg $ is of type $\textsf{B}_3$ or $\textsf{G}_2$, and $\alpha_j$ is a short root, namely, $\alpha_j = \alpha_3$ for $\textsf{B}_3$ and $\alpha_j=\alpha_1$ for $\textsf{G}_2$. In both cases $\alpha_k = \alpha_2$.  From Table 2 in page 66 of \cite{Hu}, we have $\theta = \alpha_1 + 2 \alpha_2 + 2 \alpha_3$ for $\textsf{B}_3$ and $\theta = 3 \alpha_1 + 2 \alpha_2$ for $\textsf{G}_2$. Define $$\gamma:= \theta - \alpha_2 - \alpha_j =
\left\{ \begin{array}{ll} \alpha_1 + \alpha_2 + \alpha_3 & \mbox{ for $\textsf{B}_3$} \\ 2 \alpha_1 + \alpha_2 & \mbox{ for $\textsf{G}_2$}. \end{array} \right.$$
Then it is easy to check (for example, from the list of positive roots in \cite{Ti}) that \begin{equation}\label{e.positive} \gamma, \ \gamma -\alpha_j, \ \gamma - \alpha_2 - \alpha_j \ \in \Phi^+.\end{equation}
From $$\wt(h) = 2 \theta - \alpha_2 + \beta = 2 \theta - \alpha_2 - \alpha_j= \theta + \gamma,$$ we have
$$\wt(h) - \theta + \gamma = \theta + (\gamma - \alpha_2 - \alpha_j),$$
$$\wt(h) + (- \theta + \alpha_2) + \gamma = \theta + (\gamma - \alpha_2 - \alpha_j) + \alpha_2.$$ Since neither of them can belong to $\Phi \cup \{0 \}$ by (\ref{e.positive}),
we obtain \begin{equation}\label{e.zerozero}
h(v_{-\theta}, v_{\gamma}) = h(v_{-\theta + \alpha_2}, v_{\gamma}) =0 \end{equation}
for any nonzero root vector $v_{\gamma} \in \fg_{\gamma}$.
On the other hand, Lemma \ref{l.lambda} (iii) says $$h(v_{-\theta}, v_{-\theta+ \alpha_2}) \in \C^{\times} v_{-\alpha_j}$$ for some nonzero root vector $v_{- \alpha_j} \in \fg_{-\alpha_j}$.
Combined with (\ref{e.zerozero}), this implies that $$h^{\sharp}(v_{-\theta}, v_{-\theta + \alpha_2}, v_{\gamma}) = h(v_{-\theta}, v_{-\theta+ \alpha_2}) \cdot v_{\gamma} \in \C^{\times} \, [ v_{-\alpha_j}, v_{\gamma} ] = \fg_{\gamma -\alpha_j} \setminus \{  0 \} $$ by (\ref{e.positive}),  a contradiction to $h^{\sharp} =0.$

\medskip
{\bf Case when $\delta=\theta-\alpha_k$ and $\lambda=\wt(h)=\theta$}. The argument for this case is somewhat involved. Since $\beta \neq 0$, we have $c=0$ by Lemma \ref{l.lambda} (v)  and consequently $\sigma(R)=0$ by Lemma \ref{l.lambda}.
Since $h(R)\neq 0$ by assumption, we have $f(R)\neq 0$ and \begin{equation}\label{e.ufh} u:=f(v_{-\theta}\wedge v_{-\theta+\alpha_k}) = h (v_{-\theta}\wedge v_{-\theta+\alpha_k}) \ \in\C^{\times} v_{-\theta+\alpha_k}. \end{equation}
 Let $v_{-\theta}^*$ be the unique element of $\fg^*$ satisfying $$v_{-\theta}^*(\ft) =0, \ v_{-\theta}^*(v_{-\theta}) =1 \mbox{ and }
            v_{-\theta}^*(\fg_{\alpha}) =0 \mbox{ for any root }\alpha \neq - \theta.$$
           Define $\chi \in \Hom(\wedge^2 \fg, \fg) $ by $$\chi(u, w) = v_{-\theta}^*(u) \, w - v_{-\theta}^*(w) \, u \mbox{ for all } u, w \in \fg$$ such that \begin{equation}\label{e.chi} \chi (v_{-\theta}, v_{-\theta + \alpha_k}) = v_{-\theta + \alpha_k}. \end{equation}
From (\ref{e.ufh}) and (\ref{e.chi}),  we may  replace $h$ by its suitable scalar multiple to have
$$(h-\chi)(v_{-\theta}, v_{-\theta+\alpha_k})=0.$$
Note that both $h$ and $\chi$ are ${\rm Unip}(B)$-invariant and ${\rm Unip}(B) \cdot [v_{-\theta}\wedge v_{-\theta+\alpha_k}]$ is linearly nondegenerate in $\mathbb{P}R$. It follows that
$(h-\chi)(R)=0.$
Consequently, there exist $\psi\in\End(\fg)$ and $\eta\in\fg^*$ such that
$$h-\chi=\psi \circ [\,,\,] + \eta([\,,\,]) {\rm Id}_{\fg} .$$
Furthermore, both $\psi$ and $\eta$ are ${\rm Unip}(B)$-invariant and are $\ft$-eigenvectors of weight $\theta$, since so are $h$ and $\chi$. Thus $\psi$ belongs to a $\fg$-submodule of $\End(\fg)$ isomorphic to $\fg$.  By Lemma \ref{l.decompo}, there exists $s \in\C$ such that $\psi=s \, \ad_{v_\theta}$. Since $\eta$ is a highest weight vector of $\fg^*$, it is equal to $s' \, v_{-\theta}^*$ for some $s' \in\C$.
In summary, we have
\begin{eqnarray}\label{eqn-formula-h}
h=\chi+s\, \ad_{v_{\theta}}\circ[\,,\,] + s' \, (v_{-\theta}^*\circ[\,,\,]) \otimes {\rm Id}_{\fg}
\end{eqnarray} for some  $s, s' \in\C.$
We claim \begin{equation}\label{e.hsharp}
%&& h^{\sharp}(v_{-\theta}, v_{-\theta+\alpha_k}, v_{-\alpha_k})=(-2+2s-\mu)v_{-\theta}, \\
%&& h^{\sharp}(v_{-\theta}, v_{\theta-\alpha_k}, v_{\alpha_k})=(-2-2s)v_{\theta}, \nonumber\\
  h^{\sharp}(v_{-\theta}, v_{-\alpha_k}, v_{\alpha_k})=(-2-s)t_{\alpha_k}+st_\theta,
\end{equation}
A direct calculation shows that for $x, y, z\in\fg$,
\begin{eqnarray*} \chi^{\sharp}(x, y, z)&=& 2v_{-\theta}^*(x)[y, z]+2v_{-\theta}^*(y)[z, x] +2v_{-\theta}^*(z)[x, y], \\
(\ad_{v_\theta}\circ[\,,\,])^{\sharp}(x, y, z) &=& [\ad_{v_{\theta}}([x, y]), z] +[\ad_{v_{\theta}}([y, z]), x]  \\ && +[\ad_{v_{\theta}}([z, x]), y], \\
((v_{-\theta}^*\circ[\,,\,])\otimes {\rm Id}_{\fg})^{\sharp}(x, y, z) &=& v_{-\theta}^*([x, y])z+v_{-\theta}^*([y, z])x  +v_{-\theta}^*([z, x])y.\end{eqnarray*}
Then \begin{eqnarray*}
\chi^{\sharp} (v_{-\theta}, v_{-\alpha_k}, v_{\alpha_k})
&=& 2v_{-\theta}^*(v_{-\theta})[ v_{-\alpha_k}, v_{\alpha_k}] =-2t_{\alpha_k} \\
(\ad_{v_\theta}\circ[\,,\,])^{\sharp}(v_{-\theta}, v_{-\alpha_k}, v_{\alpha_k})
&=& [\ad_{v_{\theta}}([v_{-\alpha_k}, v_{\alpha_k}]), v_{-\theta}] \\ && +[\ad_{v_{\theta}}([v_{\alpha_k}, v_{-\theta}]), v_{-\alpha_k}] \\
&= &[[v_\theta, -t_{\alpha_k}], v_{-\theta   }]+[[v_{\alpha_k}, t_\theta], v_{-\alpha_k}] \\
&=& [\langle\theta, \alpha_k\rangle v_\theta, v_{-\theta}]+[-\langle\alpha_k, \theta\rangle v_{\alpha_k}, v_{-\alpha_k}] \\
&=& t_\theta-t_{\alpha_k} \\
((v_{-\theta}^*\circ[\,,\,])\otimes {\rm Id}_{\fg})^{\sharp}(v_{-\theta}, v_{-\alpha_k}, v_{\alpha_k}) & =& 0. \end{eqnarray*}
Combing these  with \eqref{eqn-formula-h}, we obtain (\ref{e.hsharp}).
But $h^{\sharp} =0$ and (\ref{e.hsharp}) give  $s=-2$ and $s=0$, a contraction, completing the proof.
\end{proof}

We are ready to finish the proof of Proposition \ref{p.hR}.

\begin{proof}[Proof of Proposition \ref{p.hR}]
Assuming $w\neq 0$ in Lemma \ref{l.lambda},  we deduce a contradiction.
Since $h$, $v_{-\theta}$ and $v_{-\theta+\alpha_k}$ are all $\ft$-eigenvectors, we know that $w$ is a $\ft$-eigenvector of weight $\beta:=\lambda-2\theta+\alpha_k$, with $\beta\in\Sigma\cup\{0\}$ and $\wt(h)=2\theta-\alpha_k+\beta$.
We consider the following two cases separately.

\medskip
{\bf Case when $\beta\neq 0$.}  From Lemma \ref{l.alphak} and Lemma  \ref{l.exclude},   we know that the root $\beta$ must be one of the following. In each case, we choose a root $\gamma $ satisfying $\beta + \gamma \in \Phi$ and a nonzero root vector $z \in \fg_{\gamma}$ as indicated below.
\begin{itemize}
\item[$(i)$] If $\beta= \theta$, choose $z \in\fg_{\gamma}$ with $\gamma = -\alpha_k$.
\item[$(ii)$] If $\beta$ is a positive root and $\beta \neq \theta$, choose $z \in \fg_{\gamma}$ for some positive root $\gamma$ such that $\beta+\gamma\in\Phi_1\cup \Phi_2.$ In fact, if $\beta \in \Phi_1,$ we can find such $\gamma \in \Phi_1$ by Lemma \ref{l.alphak} (iii) and if $\beta \in \Phi_0$, we can find such $\gamma \in \Phi_1$  by Lemma \ref{l.alphak} (v).
\item[$(iii)$] If $\beta \in \Phi_{-1}$ and $\beta\neq-\theta+\alpha_k$, choose $z\in\fg_{\gamma}$ with $\gamma = \theta$.
\item[$(iv)$] If $\beta$ is as in Lemma \ref{l.gamma}, namely,  a negative root such that $c_k(\beta)=0$ and $-\beta$ is not simple, then we choose $z \in \fg_{\gamma}$ for the root $\gamma$ from Lemma \ref{l.gamma} satisfying $\beta + \gamma \in \Phi^-$ and $\langle \beta+ \gamma, \alpha_k \rangle \leq 0.$
\end{itemize}
Write $h = f + \sigma \, {\rm Id}_{\fg}$ as in Lemma \ref{l.lambda} and set $x=v_{-\theta}\in\fg$ and $y=v_{-\theta+\alpha_k}\in\fg$.
Then  $\sigma(x,y) =0$ by Lemma \ref{l.lambda} (v)  and $f(x,y) = h(x,y) = w \neq 0$ by our assumption.
Since $w \in \fg_{\beta}, z \in \fg_{\gamma}$ and $\beta + \gamma \in \Phi$, we have \begin{equation}\label{e.hxy} [h(x,y),z] = [f(x,y), z] = [w, z] \in \fg_{\beta + \gamma} \setminus \{0\}. \end{equation}

We claim
\begin{equation}\label{e.hyz} \wt(h)+\wt(y)+\wt(z)=\theta+\beta+\gamma\notin\Phi \cup\{0\}. \end{equation}
This can be seen from $\beta + \gamma \in \Phi^+$ in the cases (i), (ii), (iii) above.
In the case (iv), we have $c_k(\beta) = c_k (\gamma) =0$ and $c_k(\theta) =2$. Thus if $\theta + \beta + \gamma$ is a root, it must be in $\Phi_2 = \{ \theta \}$, a contradiction to $\beta + \gamma \neq 0$. On the other hand, if $\theta + \beta + \gamma =0$, then $- \beta = \theta + \gamma $ is a positive root, a contradiction to $\gamma \in \Phi^+$. This proves the claim (\ref{e.hyz}).

We claim also
\begin{equation}\label{e.hzx}
\wt(h)+\wt(z)+\wt(x)=\theta-\alpha_k+\beta+\gamma\notin\Phi\cup\{0\}.
\end{equation}
In case (i), we have $\theta - \alpha_k + \beta + \gamma = 2 \theta - 2 \alpha_k$, which cannot be zero.  If it is a root,  it must belong to $\Phi_2 = \{\theta\}$ because
$c_k( 2 \theta - 2 \alpha_k) = 2$, implying $\theta = 2 \alpha_k$, a contradiction.
In case (ii), it is clear that $\theta  - \alpha_k + \beta + \gamma \neq 0$. If $\theta - \alpha_k + \beta + \gamma \in \Phi$, it must be equal to $\theta$ from   $c_k(\theta - \alpha_k + \beta + \gamma) \geq 2$ and $\Phi_2 = \{ \theta\}.$ Then $\beta + \gamma - \alpha_k =0$, a contradiction.
In case (iii), from $c_k(\theta - \alpha_k + \beta + \gamma) = 2,$ we know that  $\theta -\alpha_k + \beta + \gamma$ is not zero. If $\theta -\alpha_k + \beta + \gamma$ is a root, then $\theta - \alpha_k + \beta  + \gamma = \theta$ from $\Phi_2 = \{ \theta\}.$ This implies
$\beta = -\theta + \alpha_k$, a contradiction to the assumption of (iii).
In case (iv), from $c_k(\theta - \alpha_k + \beta + \gamma) =1$, we know that
 $\rho:= \theta - \alpha_k + \beta + \gamma$  is not zero and if it is a root, then it is in $\Phi_1$, a positive root. Then  $\alpha_k - \beta -\gamma = \theta - \rho$ is a root. But $\langle \beta + \gamma, \alpha_k\rangle \leq 0$ implies that $\alpha_k - \beta -\gamma$ cannot be a root. This completes the proof of the claim (\ref{e.hzx}).

Since $h(y,z) \cdot x =0$ by (\ref{e.hyz}) and $h(z, x) \cdot y =0$ by (\ref{e.hzx}), we have $$ h^{\sharp}(x,y,z) = h(x,y) \cdot z \neq 0$$ from (\ref{e.hxy}), a contradiction to $h^{\sharp} =0$.

\medskip
{\bf Case when $\beta=0$.}  We claim  \begin{equation}\label{e.uf} u=f(v_{-\theta}\wedge v_{-\theta+\alpha_k}) = 0.\end{equation} Suppose the contrary, then $u$ is a nonzero element of $\ft$ from $\wt(u) = \beta =0$.  Choose two positive roots $\gamma_1$ and $\gamma_2$ different from $\alpha_k$ and root vectors $v_{\gamma_1} \in \fg_{\gamma_1}, v_{\gamma_2} \in \fg_{\gamma_2}$ such that $$[u, v_{\gamma_1}]=s_1 v_{\gamma_1} \mbox{ and } [u, v_{\gamma_2}]=s_2 v_{\gamma_2}$$ for some $s_1 \neq s_2 \in \C$.
When $x=v_{-\theta}$, $y=v_{-\theta+\alpha_k}$ and $z_i=v_{\gamma_i}, i =1, 2$, we have
\begin{equation}\label{e.hyz0} \wt(h)+\wt(y)+\wt(z_i)=\theta+\gamma_i \notin\Phi \cup\{0\}. \end{equation}
Also, \begin{equation}\label{e.hzx0}
\wt(h)+\wt(z_i)+\wt(x)=\theta-\alpha_k+\gamma_i \notin\Phi\cup\{0\},
\end{equation}
because if $\rho:= \theta - \alpha_k + \gamma_i$ is a root, then $c_k(\rho) =1$ implies that $\rho \in \Phi_1$ and  then $\alpha_k - \gamma_i = \theta -\rho \in \Phi_1,$ a contradiction.  From (\ref{e.hyz0}) and (\ref{e.hzx0}), we have $$h^{\sharp} (x, y, z_i)= h(x,y) \cdot z_i = [u, z_i] + \sigma (x,y) z_i = (s_i+\sigma(x,y))v_{\gamma_i}.$$  Then $h^{\sharp} =0$ implies $s_1=-\sigma(x,y)=s_2$, a contradiction to $s_1 \neq s_2$.
This proves (\ref{e.uf}).

Let $\gamma_1$ be a simple root different from $\alpha_k$. When $x=v_{-\theta}$, $y=v_{-\theta+\alpha_k}$ and $z=v_{\gamma_1}$,
we still have (\ref{e.hyz0}) and (\ref{e.hzx0}) with $i=1$, implying $$h^{\sharp}(x, y,z) = \sigma(x,y) z.$$ Thus $c= \sigma(v_{-\theta}, v_{-\theta + \alpha_k}) =0$.
Combining it with (\ref{e.uf}), we obtain $w=0$, which implies $h(R) =0$ by Lemma \ref{l.lambda}. This completes the proof of Proposition \ref{p.hR}.
\end{proof}

\subsection{Theorem \ref{t.Bianchi} for $\fg$ of type $\textsf{A}_{2}$}\label{ss.A2}
Throughout this subsection, the simple Lie algebra $\fg$ is of type $\textsf{A}_{2}$. Our goal is to prove

\begin{theorem}\label{t.Bianchi2}
Let $\fg$ be a simple Lie algebra of type $\textsf{A}_{2}$. Let $\ad_{\fg} \subset \End(\fg)$ be the image of the adjoint representation and denote by $\widehat{\fg} \subset \End(\fg)$  the direct product Lie algebra   $$ \widehat{\fg} :=  \ad_{\fg} \oplus \C {\rm Id}_{\fg}  \ \subset \ \End(\fg).$$ Then   $\BK(\widehat{\fg})$ is 1-dimensional, generated by the Lie bracket $$[\, ,\,] \in \Hom(\wedge^2 \fg, \ad_{\fg}) \subset \Hom(\wedge^2 \fg, \widehat{\fg}).$$
In other words, if $h \in \Hom(\wedge^2 \fg, \widehat{\fg})$  satisfies $$h^{\sharp}(x,y,z) :=  h(x,y) \cdot z + h(y,z) \cdot x + h(z,x) \cdot y  =0$$ for all $x, y, z \in \fg$, then $h = s \, [\, ,\,] \in \Hom (\wedge^2 \fg, \ad_{\fg})$ for some constant $s \in \C.$ \end{theorem}

Let us introduce some notation and recall basic properties of the Lie algebra $\fg= \fsl_3$.

\begin{notation}\label{n.sl3}
Regarding $\fg$ as the Lie algebra of traceless $3 \times 3$ matrices, let $\ft \subset \fg$ be the Cartan subalgebra of diagonal matrices (e.g. page 457 of \cite{Ya}). Denote by $\Phi \subset \ft^*$ the root system and by $\fg_{\alpha}$ the root space of a root $\alpha \in \Phi$.  Let $\lambda_i \in \ft^*, i=1,2,3,$ be the linear functional taking $i$-th entry of a diagonal matrix. We have simple roots $$\alpha_1 = \lambda_1 - \lambda_2, \ \alpha_2 = \lambda_2 - \lambda_3$$ and the highest root $\theta = \alpha_1 + \alpha_2 = \lambda_1 - \lambda_3$. The positive roots are $\Phi^+ = \{ \alpha_1, \alpha_2, \theta\}$. Denoting by $E_{ij}$ the $3\times 3$-matrix with $(i,j)$-component equal to 1 and all the other components equal to 0, we choose the following basis of $\fg$.
$$t_1 := E_{11} - E_{22}, \ t_2 := E_{22} - E_{33}, $$ $$ v_1 := E_{12}, \ v_2 = E_{23}, \ v_{\theta} = E_{13},$$  $$ v_{-1} = E_{21}, \ , v_{-2} = E_{32}, \ v_{-\theta} = -E_{31}.$$
Set $t_{\theta} := t_1 + t_2.$
 \end{notation}

The following can be checked directly.

\begin{lemma}\label{l.sl3}
 Notation \ref{n.sl3},
\begin{itemize}
\item[(i)] $v_{\theta} \in \fg_{\theta}, v_{-\theta} \in \fg_{-\theta}, v_i \in \fg_{\alpha_i}$ and $ v_{-i} \in \fg_{-\alpha_i}$ for $i=1, 2$.
\item[(ii)] For each root $\gamma\in\Phi$ and a root vector $x \in \fg_{\gamma}$,  $$
[t_1, x]=\langle \gamma, \alpha_1 \rangle x, \ [t_2, x] =\langle \gamma, \alpha_2 \rangle x, \ [t_{\theta}, x] = \langle \gamma, \theta\rangle x.$$
\item[(iii)] $t_1 = [v_1, v_{-1}], t_2 = [v_2, v_{-2}]$ and $t_{\theta} = - [ v_{\theta}, v_{-\theta}].$
\item[(iv)] We have the following relations among root vectors.
\begin{eqnarray*}
&[v_1, v_2] = v_{\theta}, & [v_{-1}, v_{-2}] = v_{-\theta} \\
& [v_{-1}, v_{\theta}]=v_2, & [v_{-2}, v_{\theta}]=-v_1, \\
& [v_1, v_{-\theta}]=v_{-2}, & [v_2, v_{-\theta}]=-v_{-1}.
\end{eqnarray*} \end{itemize}
\end{lemma}

The following can be found from Table 5 on p. 300 of \cite{OV}.

\begin{lemma}\label{l.wedgeA2}
The decomposition of $\wedge^2 \fg$ into irreducible $\fg$-modules is
$$\wedge^2 \fg = R_1 \oplus R_2 \oplus \widetilde{\fg},$$
where $R_1$ (resp. $R_2$) has $v_{-\theta} \wedge v_{-1}$ (resp. $v_{-\theta} \wedge v_{-2}$) as a lowest weight vector and  the Lie bracket $[\, , \,]: \wedge^2 \fg \to \fg$ annihilates $R_1 \oplus R_2$, sending $\widetilde{\fg}$ isomorphically to $\fg.$  \end{lemma}

\begin{lemma}\label{l.inequ}
For any nonzero dominant weight $\lambda\in \ft^*,$
there exist two positive rational numbers $r_1$ and $r_2$ such that
$$ \lambda =r_1\alpha_1+r_2\alpha_2 \mbox{ with } \frac{1}{2}\leq  \frac{r_1}{r_2}\leq 2. $$ \end{lemma}

\begin{proof}
The Lie algebra $\fg = \fsl_3$ has the fundamental weights $\pi_1 = \frac{1}{3}(2 \alpha_1 + \alpha_2)$ and $\pi_2 = \frac{1}{3}(\alpha_1 + 2 \alpha_2)$ (e.g. from Table 1 in page 294 of  \cite{OV}).  A dominant weight $\lambda$ must be of the form $ n_1 \pi_1 + n_2 \pi_2$ for some nonnegative integers $n_1, n_2$, from which we obtain $r_1 = \frac{1}{3}(2 n_1 + n_2)$ and $ r_2 = \frac{1}{3}(n_1 + 2n_2)$.
  \end{proof}

To prove Theorem \ref{t.Bianchi2}, it suffices to prove the following.

\begin{proposition}\label{p.A2}
In Theorem \ref{t.Bianchi2}, for any  irreducible  $\fg$-submodule $S \subset \BK(\widehat{\fg})$,  any highest weight vector $h \in S$ is of the form $s \, [\, , \,]$ for some $s \in \C$. \end{proposition}

We prove Proposition \ref{p.A2} in several lemmata.

\begin{lemma}\label{l.cases}
For  $h\in S$ in Proposition \ref{p.A2}, write
\begin{eqnarray*}
u_0:=h(v_{-1}, v_{-2}), & u_1:=h(v_{-\theta}, v_{-1}), & u_2:=h(v_{-\theta}, v_{-2}).
\end{eqnarray*}
Then the following holds.
\begin{itemize}
\item[(i)] All of $u_0, u_1$ and $u_2$ are $\ft$-eigenvectors.
\item[(ii)] $h=0$ if and only if $u_0=u_1=u_2=0$.
\item[(iii)] Let $\beta_0,\beta_1,\beta_2$ and $\lambda$ be the weights of $u_0, u_1, u_2$ and $h$ respectively. If $h\neq 0$, then we have the following seven possibilities of the 4-tuple  $(\lambda, \beta_0, \beta_1, \beta_2)$. Here we write $\infty$ for the weight when the corresponding weight vector is zero (thus $u_1 = u_2 =0$ in the case 1, etc.).
\begin{itemize}
\item[1.] $(0, -\alpha_1-\alpha_2, \infty, \infty)$;
\item[2.] $(2\alpha_1+\alpha_2, \alpha_1, 0, \infty)$;
\item[3.] $(\alpha_1+2\alpha_2, \alpha_2, \infty, 0)$;
\item[4.] $(2\alpha_1+2\alpha_2, \alpha_1+\alpha_2, \alpha_2, \alpha_1)$;
\item[5.] $(3\alpha_1+2\alpha_2, \infty, \alpha_1+\alpha_2, \infty)$;
\item[6.] $(2\alpha_1+3\alpha_2, \infty, \infty, \alpha_1+\alpha_2)$;
\item[7.] $(\alpha_1+\alpha_2, 0, -\alpha_1, -\alpha_2)$.
\end{itemize}
\end{itemize}
\end{lemma}

\begin{proof}
(i) is straightforward and  (ii) is a direct consequence of Lemma \ref{l.wedgeA2}.
To check (iii), note that  at least one of $u_0$, $u_1$ and $u_2$ is nonzero by (ii).
If $u_1\neq 0$, then it is a $\ft$-eigenvector of $\fg\oplus \C {\rm Id}_{\fg},$ implying $\beta_1\in\Phi \cup\{0\}$. Since $h$ is a highest weight vector, its weight $\lambda$ satisfies the condition from Lemma \ref{l.inequ}.  From these two facts, we see that  $(\lambda, \beta_0, \beta_1, \beta_2)$ should be one of the cases  2, 4, 5 and 7. Arguing in a  similar way  for the cases $u_0\neq 0$ and $u_2\neq 0$, we see that one of the cases 1-7 occurs.
\end{proof}

To prove Proposition \ref{p.A2}, we check $h \in \C \, [ \, , \,]$ case-by-case for each  of the seven possibilities in Lemma \ref{l.cases} (iii). In the computation below, we use the relations from Lemma \ref{l.sl3}.

\begin{lemma}\label{l.case1}
In case 1 of Lemma \ref{l.cases}, we have $h\in\C[\,,\,]$.
\end{lemma}

\begin{proof}
In this case, we have
$ h(R_1)=h(R_2)=0$ and $ u_0:=h(v_{-1}, v_{-2})=sv_{-\theta}$ for some $ s\in\C. $
Define $\psi:=h-s[\,,\,].$ Then $\psi( v, v') =0$ for all negative root vectors $v, v' \in \fg_{-\alpha_1} \oplus \fg_{- \alpha_2} \oplus \fg_{-\theta}$. Since $h$ is a highest weight vector of weight $\lambda =0$, it is invariant under the action of the unipotent radical ${\rm Unip}(B)$ of the Borel subgroup $B$. Thus $\psi$ is also invariant under ${\rm Unip}(B)$. It follows that $\psi(v, v') =0$ for all $v, v' \in \fg$. \end{proof}

\begin{lemma}\label{l.case2}
In cases 2 and 3 of Lemma \ref{l.cases}, we have $h=0$.
\end{lemma}

\begin{proof}
By the symmetry of the Dynkin diagram of $\textsf{A}_2$, it suffices to prove it for case 2. For case 2,  we have  $u_0:=h(v_{-1}, v_{-2})=s_0 v_1$ for some $ s_0\in\C$, $u_2:=h(v_{-\theta}, v_{-2})=0$ and
\begin{equation}\label{eqn-case2-tg1}
 u_1:=h(v_{-\theta}, v_{-1})=\zeta+s_1
\end{equation}
for some $ \zeta\in\ft $ and some $s_1\in\C. $
Since $v_2\cdot h=0$ and $v_2\cdot(v_{-\theta}\wedge v_{-1})=0$, we have $v_2\cdot u_1=0$, which implies that $\zeta\in\C (2t_1+t_2)$. Then we can rewrite \eqref{eqn-case2-tg1} as
$$ u_1:=h(v_{-\theta}, v_{-1})=2s_2t_1+s_2t_2+s_1
$$ for some $ s_1, s_2\in\C.$
A direct calculation shows that
\begin{eqnarray}
&&\hspace{0.5cm}  h^{\sharp} (v_{-\theta}, v_{-1}, v_{-2}) \label{eqn-case2-dh1} \\
&& = h(v_{-\theta}, v_{-1})\cdot v_{-2}+h(v_{-1}, v_{-2})\cdot v_{-\theta}+ h(v_{-2}, v_{-\theta})\cdot v_{-1} \nonumber\\ && =(2s_2t_1+s_2t_2+s_1)\cdot v_{-2}+[s_0v_1, v_{-\theta}]+0\cdot v_{-1} \nonumber\\
&& =(s_0+s_1)v_{-2}, \nonumber\\
&& \nonumber\\
&&\hspace{0.5cm}   h^{\sharp}(v_{-\theta}, v_1, v_{-1}) \label{eqn-case2-dh2} \\
&& = h(v_{-\theta}, v_1)\cdot v_{-1}+h(v_1, v_{-1})\cdot v_{-\theta}+ h(v_{-1}, v_{-\theta})\cdot v_1 \nonumber\\
&& = 0\cdot v_{-1}+0\cdot v_{-\theta}-(2s_2t_1+s_2t_2+s_1)\cdot v_1 \nonumber\\
&& = (s_1+3 s_2)v_1,\nonumber\\
&& \nonumber\\
&&\hspace{0.5cm}  h^{\sharp}(v_2, v_{-1}, v_{-2}) \label{eqn-case2-dh3} \\
&& = h(v_2, v_{-1})\cdot v_{-2}+h(v_{-1}, v_{-2})\cdot v_2+ h(v_{-2}, v_2)\cdot v_{-1} \nonumber\\
&& = 0\cdot v_{-2}+[s_0v_1, v_2]+ 0\cdot v_{-1} \nonumber\\
&& = s_0v_\theta. \nonumber
\end{eqnarray}
Here, the vanishing of $h(v_{-2}, v_{-\theta}), h(v_{-\theta}, v_1), h(v_1, v_{-1}), h(v_2, v_{-1})$ and $ h(v_{-2}, v_2),$ which are $\ft$-eigenvectors in $\fg$,  comes from the fact that their weights cannot be  roots.
From the vanishing of  \eqref{eqn-case2-dh1} \eqref{eqn-case2-dh2} and \eqref{eqn-case2-dh3}, we have $s_0=s_1=s_2=0$. It follows that $u_0=u_1=u_2=0$, implying $h=0$ by Lemma \ref{l.cases}.
\end{proof}

\begin{lemma}\label{l.case4}
In case 4 of Lemma \ref{l.cases}, we have $h=0$.
\end{lemma}

\begin{proof}
In this case, there exist $s_0, s_1, s_2\in\C$ such that
\begin{eqnarray*}
u_0 &:= & h(v_{-1}, v_{-2})=s_0v_\theta, \\ u_1 &:=& h(v_{-\theta}, v_{-1})=s_1v_2, \\  u_2 &:= & h(v_{-\theta}, v_{-2})=s_2v_1.
\end{eqnarray*}
A direct calculation shows that
\begin{eqnarray}
&&\hspace{0.5cm}  h^{\sharp}(v_{-\theta}, v_{-1}, v_{-2}) \label{eqn-case4-dh1} \\
&& = h(v_{-\theta}, v_{-1})\cdot v_{-2}+h(v_{-1}, v_{-2})\cdot v_{-\theta}+ h(v_{-2}, v_{-\theta})\cdot v_{-1} \nonumber\\
&& =[s_1v_2, v_{-2}]+[s_0v_\theta, v_{-\theta}]+[-s_2v_1, v_{-1}] \nonumber\\
&& =s_1t_2-s_0t_\theta - s_2t_1, \nonumber\\
&& =-(s_0+s_2)t_1+(s_1-s_0)t_2, \nonumber\\
&& \nonumber\\
&&\hspace{0.5cm}  h^{\sharp}(v_{-\theta}, v_1, v_{-1}) \label{eqn-case4-dh2} \\
&& = h(v_{-\theta}, v_1)\cdot v_{-1}+h(v_1, v_{-1})\cdot v_{-\theta}+ h(v_{-1}, v_{-\theta})\cdot v_1 \nonumber\\
&& = 0\cdot v_{-1}+0\cdot v_{-\theta}+[-s_1v_2, v_1] \nonumber\\
&& = s_1v_\theta \nonumber
\end{eqnarray}
Here, the vanishing of $h(v_{-\theta}, v_1)$ and $h(v_1, v_{-1})$ is from the fact that their weights cannot be roots.
From the vanishing of \eqref{eqn-case4-dh1} and \eqref{eqn-case4-dh2}, we have $s_0=s_1=s_2=0$. It follows that $u_0=u_1=u_2=0$, implying $h=0$ by Lemma \ref{l.cases}.
\end{proof}

\begin{lemma}\label{lem-case5}
In case 5 and case 6 of Lemma \ref{l.cases}, we have $h=0$.
\end{lemma}

\begin{proof}
By the symmetry of the Dynkin diagram of $\textsf{A}_2$, it suffices to prove it for case 5. For case 5, we have
\begin{eqnarray*}
u_0:=h(v_{-1}, v_{-2})=0, & u_1:=h(v_{-\theta}, v_{-1})=sv_\theta, & u_2:=h(v_{-\theta}, v_{-2})=0
\end{eqnarray*}  for some  $s\in\C$.
A direct calculation shows that
\begin{eqnarray*}
&&\hspace{0.5cm} h^{\sharp}(v_{-\theta}, v_{-1}, v_{-2})  \\
&& = h(v_{-\theta}, v_{-1})\cdot v_{-2}+h(v_{-1}, v_{-2})\cdot v_{-\theta}+ h(v_{-2}, v_{-\theta})\cdot v_{-1}\\
&& =[sv_\theta, v_{-2}]+0\cdot v_{-\theta}+0\cdot v_{-1} \\
&& =sv_1.
\end{eqnarray*}
Here, the vanishing of $h(v_{-1}, v_{-2})$ and $h(v_{-2}, v_{-\theta})$ is from the fact that their weights cannot be roots.
From $h^{\sharp} =0$,  we obtain $s=0$. It follows that $u_0=u_1=u_2=0$, implying $h=0$ by Lemma \ref{l.cases}.
\end{proof}

\begin{lemma}\label{l.case7}
In case 7 of Lemma \ref{l.cases}, we have $h=0$.
\end{lemma}

\begin{proof}
In this case there are constants $\xi_0, \xi_1, \xi_2, s_1, s_2\in\C$ such that
\begin{eqnarray}
&& h(v_{-\theta}, v_{-1})=s_1v_{-1}\in\fg_{-\alpha_1}, \label{eqn2}\\
&& h(v_{-\theta}, v_{-2})=s_2v_{-2}\in\fg_{-\alpha_2}, \label{eqn3}\\
&& h(v_{-1}, v_{-2})=\xi_1t_1+\xi_2t_2+\xi_0 \, {\rm Id}_{\fg} \in\ft\oplus\C {\rm Id}_{\fg}. \label{eqn4}
\end{eqnarray}
The proof of Lemma \ref{l.case7} is in two steps.
In Step 1, we determine the linear map $h$ completely. In Step 2, we use the results of Step 1 to compute  $h^{\sharp}$ to deduce $\xi_i = s_j =0$, from which $h=0$ follows.

\medskip

{\bf Step 1.} Applying the action of $v_1\in\fg$ on both sides of \eqref{eqn2}, we have
\begin{eqnarray}
-h(v_{-1}, v_{-2})+h(v_{-\theta}, t_1)=s_1t_1.\label{eqn5}
\end{eqnarray}
Applying the action of $v_2\in\fg$ on both sides of \eqref{eqn3}, we have
\begin{eqnarray}
-h(v_{-1}, v_{-2})+h(v_{-\theta}, t_2)=s_2t_2.\label{eqn6}
\end{eqnarray}
By \eqref{eqn4} and \eqref{eqn5}, we have
\begin{eqnarray}
h(v_{-\theta}, t_1)=(\xi_1+s_1)t_1+\xi_2t_2+\xi_0 \, {\rm Id}_{\fg}. \label{eqn7}
\end{eqnarray}
By \eqref{eqn4} and \eqref{eqn6}, we have
\begin{eqnarray}
h(v_{-\theta}, t_2)=\xi_1t_1+(\xi_2+s_2)t_2+\xi_0 \, {\rm Id}_{\fg}. \label{eqn8}
\end{eqnarray}
Applying the action of $v_1\in\fg$ on both sides of \eqref{eqn7}, we have
\begin{eqnarray}
h(v_{-2}, t_1)-2h(v_{-\theta}, v_1)=(-2s_1-2\xi_1+\xi_2)v_1.\label{eqn9}
\end{eqnarray}
Applying the action of $v_2\in\fg$ on both sides of \eqref{eqn8}, we have
\begin{eqnarray}
-h(v_{-1}, t_2)-2h(v_{-\theta}, v_2)=(-2s_2+\xi_1-2\xi_2)v_2.\label{eqn10}
\end{eqnarray}
Applying the action of $v_2\in\fg$ on both sides of \eqref{eqn7}, we have
\begin{eqnarray}
-h(v_{-1}, t_1)+h(v_{-\theta}, v_2)=(s_1+\xi_1-2\xi_2)v_2.\label{eqn11}
\end{eqnarray}
Applying the action of $v_1\in\fg$ on both sides of \eqref{eqn8}, we have
\begin{eqnarray}
h(v_{-2}, t_2)+h(v_{-\theta}, v_1)=(s_2-2\xi_1+\xi_2)v_1.\label{eqn12}
\end{eqnarray}
Applying the action of $v_1\in\fg$ on both sides of \eqref{eqn4}, we have
\begin{eqnarray}
h(v_{-2}, t_1)=(2\xi_1-\xi_2)v_1.\label{eqn13}
\end{eqnarray}
Applying the action of $v_2\in\fg$ on both sides of \eqref{eqn4}, we have
\begin{eqnarray}
h(v_{-1}, t_2)=(\xi_1-2\xi_2)v_2.\label{eqn14}
\end{eqnarray}
By \eqref{eqn9} and \eqref{eqn13}, we have
\begin{eqnarray}
h(v_{-\theta}, v_1)=(s_1+2\xi_1-\xi_2)v_1.\label{eqn15}
\end{eqnarray}
By \eqref{eqn10} and \eqref{eqn14}, we have
\begin{eqnarray}
h(v_{-\theta}, v_2)=(s_2-\xi_1+2\xi_2)v_2.\label{eqn16}
\end{eqnarray}
By \eqref{eqn12} and \eqref{eqn15}, we have
\begin{eqnarray}
h(v_{-2}, t_2)=(-s_1+s_2-4\xi_1+2\xi_2)v_1.\label{eqn17}
\end{eqnarray}
By \eqref{eqn11} and \eqref{eqn16}, we have
\begin{eqnarray}
h(v_{-1}, t_1)=(-s_1+s_2-2\xi_1+4\xi_2)v_2.\label{eqn18}
\end{eqnarray}
Applying the action of $v_2\in\fg$ on both sides of \eqref{eqn13}, we have
\begin{eqnarray}
-h(t_1, t_2)-h(v_2, v_{-2})=-(2\xi_1-\xi_2)v_\theta.\label{eqn21}
\end{eqnarray}\\
Applying the action of $v_1\in\fg$ on both sides of \eqref{eqn14}, we have
\begin{eqnarray}
h(t_1, t_2)-h(v_1, v_{-1})=(\xi_1-2\xi_2)v_\theta.\label{eqn22}
\end{eqnarray}
Applying the action of $v_2\in\fg$ on both sides of \eqref{eqn15}, we have
\begin{eqnarray}
h(v_1, v_{-1})+h(v_\theta, v_{-\theta})=(-s_1-2\xi_1+\xi_2)v_\theta.\label{eqn23}
\end{eqnarray}
Applying the action of $v_1\in\fg$ on both sides of \eqref{eqn16}, we have
\begin{eqnarray}
-h(v_2, v_{-2})-h(v_\theta, v_{-\theta})=(s_2-\xi_1+2\xi_2)v_\theta.\label{eqn24}
\end{eqnarray}
By \eqref{eqn21}-\eqref{eqn24}, we have
\begin{eqnarray}
&& h(v_1, v_{-1})=(\frac{-s_1+s_2}{2}-\xi_1+2\xi_2)v_\theta, \label{eqn27} \\
&& h(v_2, v_{-2})=(\frac{s_1-s_2}{2}+2\xi_1-\xi_2)v_\theta, \label{eqn28} \\
&& h(t_1, t_2)=\frac{-s_1+s_2}{2}v_\theta, \label{eqn29} \\
&& h(v_\theta, v_{-\theta})=-(\frac{s_1+s_2}{2}+\xi_1+\xi_2)v_\theta. \label{eqn30}
\end{eqnarray}

\medskip
{\bf Step 2.} A direct calculation shows that
\begin{eqnarray}
&&h^{\sharp}(v_{-\theta}, t_1, t_2) \label{eqn36} \\
&=&  h(v_{-\theta}, t_1)\cdot t_2+h(t_1, t_2)\cdot v_{-\theta}+ h(t_2, v_{-\theta})\cdot t_1 \nonumber\\ &=& ((s_1+\xi_1)t_1+\xi_2t_2+\xi_0 \, {\rm Id}_{\fg})\cdot t_2+[\frac{-s_1+s_2}{2}v_\theta, v_{-\theta}] \nonumber \\ & & -(\xi_1t_1+(s_2+\xi_2)t_2+\xi_0 \, {\rm Id}_{\fg})\cdot t_1 \nonumber\\
&=& \xi_0t_2+\frac{s_1-s_2}{2}t_\theta-\xi_0t_1 , \nonumber\\
&=&(\frac{s_1-s_2}{2}+\xi_0)t_\theta-2\xi_0t_1 , \nonumber\\
&& \nonumber\\
&&h^{\sharp}(v_{-1}, v_{-2}, v_\theta) \label{eqn42} \\
&=& h(v_{-1}, v_{-2})\cdot v_\theta+h(v_{-2}, v_\theta)\cdot v_{-1}+ h(v_\theta, v_{-1})\cdot v_{-2} \nonumber\\ &=& (\xi_1t_1+\xi_2t_2+\xi_0\, {\rm Id}_{\fg})\cdot v_\theta+0\cdot v_{-1}+0\cdot v_{-2}\nonumber\\
&=& (\xi_1+\xi_2+\xi_0)v_\theta, \nonumber\\
&& \nonumber\\
&& h^{\sharp}(v_{-\theta}, v_1, v_{-1}) \label{eqn43} \\
&=&  h(v_{-\theta}, v_1)\cdot v_{-1}+h(v_1, v_{-1})\cdot v_{-\theta}+ h(v_{-1}, v_{-\theta})\cdot v_1 \nonumber\\
&=& [(s_1+2\xi_1-\xi_2)v_1, v_{-1}] \nonumber \\ & & +[(\frac{-s_1+s_2}{2}-\xi_1+2\xi_2)v_\theta, v_{-\theta}]+[-s_1v_{-1}, v_1] \nonumber\\
&=& (s_1+2\xi_1-\xi_2)t_1+(\frac{s_1-s_2}{2}+\xi_1-2\xi_2)t_\theta+s_1t_1 , \nonumber\\
&=& (2s_1+2\xi_1-\xi_2)t_1+(\frac{s_1-s_2}{2}+\xi_1-2\xi_2)t_\theta. \nonumber
\end{eqnarray}
From the vanishing of \eqref{eqn36} \eqref{eqn42} and \eqref{eqn43}, we have
\begin{eqnarray*}
&& \frac{s_1-s_2}{2}+\xi_0=0, \\
&& \xi_0=0, \\
&& \xi_1+\xi_2+\xi_0=0, \\
&& 2s_1+2\xi_1-\xi_2=0, \\
&& \frac{s_1-s_2}{2}+\xi_1-2\xi_2=0.
\end{eqnarray*}
It follows that $s_1=s_2=\xi_1=\xi_2=\xi_0=0$, implying $h=0$ by Lemma \ref{l.cases}.
\end{proof}

\section{Characteristic conic connections arising from minimal rational curves }\label{s.vmrt}

Many interesting examples of cone structures with characteristic conic connections arise from algebraic geometry as VMRT-structures of minimal rational curves. We have already used VMRT implicitly for adjoint varieties in Section \ref{s.adjoint}.  We  recall quickly some basic terminology and relevant results to explain the implication of Corollary \ref{c.adjoint} and Theorem \ref{t.Bianchi} for
minimal rational curves. We recommend \cite{Hw01} for a detailed introduction to VMRT.

\begin{definition}\label{d.vmrt}
Let ${\rm RatCurves}(X)$ be the space of rational curves on a smooth projective variety $X$. An irreducible component $\sK$ of ${\rm RatCurves}(X)$ is a {\em family of minimal rational curves } on $X$, if for a general point $x \in X$, the subscheme $\sK_x \subset \sK$ consisting of members of $\sK$ passing through $x$ is nonempty and projective. Fix such a $\sK$. Let $\sC^+ \subset \BP TX$ be the closure of the union of tangent vectors to members of $\sK$. By the irreducibility of $\sK$, there exists a unique irreducible component  $\sC \subset \BP TX$ of $\sC^+$ that is dominant over $X$. We call $\sC$ a {\em VMRT-structure} on $X$ (VMRT is the abbreviation of Variety of Minimal Rational Tangents). The fiber $\sC_x \subset \BP T_x X$ at a point $x \in X$ is the {\em VMRT} at $x$ of the family $\sK$.
\end{definition}

In general, the VMRT $\sC_x$ may not be smooth for a general $x \in X$. But it is smooth in many cases. For example, when the members of $\sK$ are lines under a projective embedding $X \subset \BP^N$, the VMRT at a general point is smooth. The following is from Proposition 8 of \cite{HM04} (see also   Section 3 of \cite{HN} for an explanation for differential geometers).

\begin{proposition}\label{p.conic}
Assume that the VMRT $\sC_x \subset \BP T_x X$ is smooth for a general $x \in X$, in other words,  on a Zariski-open subset $M \subset X$, the restriction $\sC|_M \subset \BP TM$ is a cone structure in the sense of Definition \ref{d.cone}. Then the tangent directions to members of $\sK$ determine a characteristic conic connection on $\sC|_M$.
\end{proposition}

Proposition \ref{p.conic} provides many interesting isotrivial cone structures equipped with characteristic connections via VMRT. The simplest examples are the following.

\begin{example}\label{ex.flat} Let   $\BP^{n-1} \subset \BP^{n}$ be a hyperplane in the $n$-dimensional projective space and fix a submanifold $Z \subset \BP^{n-1}$. Let $X_Z$ be the blowup of $\BP^n$ along $Z$ and let $\sK_Z$ be the family of minimal rational curves on $X$ whose general members are proper transformations of lines on $\BP^n$ intersecting $Z$ (see Example 1.7 of \cite{Hw10}). Then its VMRT structure restricted to the open subset of $X_Z$ corresponding to $\BP^n \setminus \BP^{n-1}$ is a  locally flat $Z$-isotrivial cone structure. \end{example}

The following result from \cite{HM01}, called the Cartan-Fubini type extension theorem, shows that VMRT-structures determine the geometry of the projective variety $X$ to some extent.

\begin{theorem}\label{t.CF}
Let $\sK$ (resp. $\widetilde{\sK}$) be a family of minimal rational curves on a smooth projective variety $X$ (resp. $\widetilde{X}$). Assume that the VMRT $\sC_x \subset \BP T_x X$ is smooth, irreducible and not a linear subspace of $\BP T_x X$.
Suppose that  there exists a biholomorphic map $\varphi: U \to \widetilde{U}$ between some connected open subsets $U \subset X$ and $\widetilde{U} \subset \widetilde{X}$ such that its differential ${\rm d} \varphi: \BP TU \to \BP T\widetilde{U}$ sends $\sC|_U$ biholomorphically to $\widetilde{\sC}|_{\widetilde{U}}.$ Then there are \begin{itemize} \item a  member $C \subset X$ of $\sK$ and a neighborhood $C \subset O \subset X$;
\item a  member $\widetilde{C} \subset \widetilde{X}$ of $\widetilde{\sK}$ and a neighborhood $\widetilde{C} \subset \widetilde{O} \subset \widetilde{X}$;
    \item a biholomorphic map $f: O \to \widetilde{O}$ with $f(C) = \widetilde{C}. $ \end{itemize} If furthermore, the second Betti numbers of $X$ and $\widetilde{X}$ are 1, then we can choose $f$ as a biholomorphic map from $X$ to $\widetilde{X}$ such that $\varphi= f|_U$. \end{theorem}

  Especially,  when  $\sC|_M$ is $Z$-isotrivial for a highest weight variety $Z \subset \BP V$ of certain types, we have some rigidity results.
     The earliest result in this direction is the following  from  Main Theorem in Section 2 of \cite{Mok}.

        \begin{theorem}\label{t.Mok}
        Let $Z \subset \BP V$ be one of the highest weight varieties listed in Proposition \ref{p.KN}. For a VMRT-structure $\sC \subset \BP TX$ on a smooth projective variety $X$, if its restriction $\sC|_M$ to some open subset $M \subset X$  is $Z$-isotrivial, then the cone structure $\sC|_M$ is locally flat. In particular, a general member of the family of minimal rational curve has a neighborhood biholomorphic to a neighborhood of a general member of the family of minimal rational curve $\sK_Z$ on $X_Z$ in Example \ref{ex.flat}. \end{theorem}

        In the setting of Theorem \ref{t.Mok}, the associated $G$-structure on $M$ is a parabolic geometry in the sense of  \cite{CS}. The proof in \cite{Mok} uses this parabolic geometry combined with the geometry of rational curves. In \cite{HN}, an alternative proof is given which uses essentially only parabolic geometry, replacing the geometry of rational curves by suitable conditions on the curvature of the Cartan connection of the parabolic geometry.

        Another result of this kind is the following from  Theorem 1.5 of \cite{HL}.

        \begin{theorem}\label{t.HL}
        Let $Z \subset \BP V$ be one of the highest weight varieties listed in Proposition \ref{p.Legendre}, excluding the twisted cubic curve. For a VMRT-structure $\sC \subset \BP TX$ on a smooth projective variety $X$, if its restriction $\sC|_M$ to some open subset $M \subset X$ is $Z$-isotrivial, then the cone structure $\sC|_M$ is locally flat. In particular, a general member of the family of minimal rational curve has a neighborhood biholomorphic to a neighborhood of a general member of the family of minimal rational curve $\sK_Z$ on $X_Z$ in Example \ref{ex.flat}. \end{theorem}

        The proof of Theorem \ref{t.HL} is rather indirect. It uses a certain induced  differential-geometric structure on $\sK$ instead of the cone structure $\sC|_M$ itself, a kind of Penrose transformation of $\sC|_M$. A Cartan connection is constructed for this structure on $\sK$ and the local flatness is checked for this connection.

Let us look at the case when the highest weight variety is the adjoint variety.
        In this case, we have more interesting examples than Example \ref{ex.flat}, as explained below.

        \begin{example}\label{ex.BF}
        Let $\fg$ be a simple Lie algebra of type different from $\textsf{A}_{\ell}$ and let $Y \subset \BP \fg$ be its adjoint variety. Let $G$ be the adjoint group of $\fg$.  There is a smooth projective variety $X_{\fg}$, called the wonderful compactification of $G$, which is equipped with an action of $G \times G$ such that there is an open orbit $M_{\fg} \subset X_{\fg}$ isomorphic to $G$ and the action of $G\times G$ on $M_{\fg}$ is given by the left and the right translations of $G$.   By the result of Brion and Fu in \cite{BF},  there exists a  unique family $\sK_{\fg}$ of minimal rational curves on $X_{\fg}$ such that its VMRT-structure restricted to an open subset of $M_{\fg}$ is a $Y$-isotrivial cone structure. It is not locally flat if $\fg$ is different from $\textsf{C}_{\ell}$. \end{example}

\begin{example}\label{ex.A}
Let ${\rm Gr}(k, \C^{2k}) \subset \BP (\wedge^k \C^{2k}), k \geq 3,$ be the Pl\"ucker embedding of the Grassmannian of $k$-dimensional subspaces in a $2k$-dimensional vector space.
Let $X^H = {\rm Gr}(k, \C^{2k})\cap H$ be the intersection with a hyperplane $H \subset \BP (\wedge^k \C^{2k}).$  For a general $H$, the family of lines on $X^H$ determines a family of minimal rational curves on $X^H$ such that  the VMRT $\sC_x \subset \BP T_x X^H$ at a general point $x \in X^H$ is isomorphic to a general hyperplane section of the VMRT associated to the family of  lines on the Grassmannian ${\rm Gr}(k, \C^{2k})$. The VMRT of the Grassmannian at a general point is isomorphic to the Segre variety $\BP^{k-1} \times \BP^{k-1} \subset \BP^{k^2-1}$ and its general hyperplane section is isomorphic to the adjoint variety of $\fsl_k$, as discussed in Example \ref{ex.sl}.
In  Theorem 1.2 of \cite{BFM}, it was proved that $\Aut(\widehat{X}^H)$ has an open orbit in $X^H$ for a general $H$ if and only if $k =3$.

When $k=3$, the adjoint variety $Y \subset \BP \fg$ is of $\textsf{A}_2$-type. Thus the VMRT of $X^H$ gives  a  $Y$-isotrivial cone structure on an open subset of $X^H$ which has a characteristic conic connection, hence is locally symmetric by Corollary \ref{c.adjoint}. It is not locally flat.

 We claim that  when $k \geq 4$, the isotrivial cone structure induced by the VMRT-structure on any open subset $M \subset X^H$ cannot be locally symmetric for a general $H$. Suppose that it is locally symmetric. Then  $\Aut(\widehat{X}^H)$ has an open orbit in $X^H$ for a general $H$ by Corollary \ref{c.homogeneous} and Theorem \ref{t.CF} because   the second Betti number of $X^H$  is 1. This is  a contradiction to Theorem 1.2 of  \cite{BFM}.
From the claim, the hyperplane section $X^H$ for $k \geq 4$ shows that the conclusion of  Corollary \ref{c.adjoint}  fails when $\fg$ is of type $\textsf{A}_{\ell \geq 3}$. \end{example}

    Theorem \ref{t.Bianchi} and Corollary \ref{c.adjoint}  have the following consequence.

        \begin{theorem}\label{t.vmrt}
        Let $Y \subset \BP \fg$ be an adjoint variety of type different from $\textsf{A}_{\ell \neq 2}$ or $\textsf{C}_{\ell}$. Let $X$ be a smooth projective variety of dimension equal to $\dim \fg$ with a family $\sK$ of minimal rational curves such that the   restriction $\sC|_M$ of the VMRT-structure  to some open subset $M \subset X$ is $Y$-isotrivial. \begin{itemize}  \item[(i)] Then the cone structure $\sC|_M \subset \BP TM$ is locally symmetric, hence locally homogeneous. \item[(ii)] A neighborhood of a general member of $\sK$ in $X$ is biholomorphic to a neighborhood of a general member of the family of minimal rational curves given in Examples \ref{ex.flat}, \ref{ex.BF} or \ref{ex.A}. \item[(iii)]
             If furthermore, the second Betti number of $X$ is 1, then $X$ is quasi-homogeneous, namely, the automorphism group of $X$ acts on $X$ with an open orbit.\end{itemize}
              \end{theorem}

\begin{proof}
 (i) follows from
Corollaries \ref{c.homogeneous} and \ref{c.adjoint}, combined with Proposition \ref{p.conic}. (ii) follows (i) combined with Theorem \ref{t.Bianchi}, Corollary \ref{c.equivcone} and Theorem \ref{t.CF}.
(iii) follows from (i) by Corollary \ref{c.homogeneous} and Theorem \ref{t.CF}.
\end{proof}

        \begin{remark}\label{r.C} It is  interesting  that Example \ref{ex.BF} works also for $\fg$ of type $\textsf{C}_{\ell}$. In that case, the adjoint variety $Y \subset \BP \fg$ is isomorphic to the case (III) of Proposition \ref{p.KN} and the VMRT-structure $\sC \subset \BP TX$ on the wonderful compactification $X$  gives a $Y$-isotrivial cone structure on an open subset $G \cong M \subset X$. Then $\sC$ is locally flat by Theorem \ref{t.Mok}. But the $G$-structure of $\sC|_M$ has a locally symmetric principal connection which is not locally flat. In other words, there are two distinct locally symmetric principal connections, one locally flat and the other not locally flat. This is possible because the Spencer homomorphism is not injective for the cases in Proposition \ref{p.KN}, hence there exist more than one torsion-free principal connections. \end{remark}

\textbf{Declarations of interest:}  none.

\bigskip
Jun-Muk Hwang (jmhwang@ibs.re.kr)
Center for Complex Geometry,
Institute for Basic Science (IBS),
Daejeon 34126, Republic of Korea

\medskip
Qifeng Li (qifengli@sdu.edu.cn)
School of Mathematics,
Shandong University,
Jinan 250100, China
\end{document}